\documentclass{article}
\usepackage{array}
\usepackage{amsmath}
\usepackage{amssymb}
\usepackage{bbm} 
\usepackage{amsfonts}
\usepackage{graphicx}
\usepackage{url}
\usepackage{maplestd2e}

\DefineParaStyle{Maple Heading 1}
\DefineParaStyle{Maple Text Output}
\DefineParaStyle{Maple Dash Item}
\DefineParaStyle{Maple Bullet Item}
\DefineParaStyle{Maple Normal}
\DefineParaStyle{Maple Heading 4}
\DefineParaStyle{Maple Heading 3}
\DefineParaStyle{Maple Heading 2}
\DefineParaStyle{Maple Warning}
\DefineParaStyle{Maple Title}
\DefineParaStyle{Maple Error}
\DefineCharStyle{Maple Hyperlink}
\DefineCharStyle{Maple 2D Math}
\DefineCharStyle{Maple Maple Input}
\DefineCharStyle{Maple 2D Output}
\DefineCharStyle{Maple 2D Input}
\usepackage[top=.9in,bottom=.9in]{geometry}
 \numberwithin{equation}{section}
 
\DefineParaStyle{Maple Output}
\DefineCharStyle{2D Comment}
\DefineCharStyle{2D Math}
\DefineCharStyle{2D Output}
\begin{document}
\pagestyle{plain}
\numberwithin{equation}{section}
\begin{flushleft}

\centerline{\bf \large{Master Theorems for a Family of Integrals}} \vskip .2in

\centerline{ M.L. Glasser\footnote{laryg@clarkson.edu}}\vskip .2in

\centerline{Department of Physics, Clarkson University}
 
\centerline{Potsdam, NY 13699-5820}
\centerline{}
\centerline{and}\vskip .2in
\centerline{ Michael Milgram\footnote{mike@geometrics-unlimited.com}}

\centerline{Consulting Physicist, Geometrics Unlimited, Ltd.}
\centerline{Box 1484, Deep River, Ont. Canada. K0J 1P0}
\centerline{}
\centerline{May 8, 2014}
\centerline{}
\vskip .5in
\centerline{\bf Abstract}\vskip .1in

A family of general Master theorems for analytic integration over the real (or imaginary) axis with various reciprocal hyperbolic (trig) kernels ($ 1/\sinh$ \it{and/or} \rm $ 1/\cosh$) with varying arguments is developed. Several examples involving closed-form integrations which do not appear to exist in the standard tables are given in detail and general results are listed in an Appendix. In addition, it is shown how to convert the special case of an infinite series involving ratios of Gamma functions into a finite series. 
 

\section{Introduction}

Recently, one of us has obtained \cite{Remarkable} a result that provides for the analytic integration of an (almost) arbitrary function, coupled with a specific kernel function, over the infinite real axis. The question naturally arose whether other kernel functions exist that yield similar formulae that allow for the integration of (almost) arbitrary functions, and, if so, is it possible to find a general procedure to discover such master kernels and determine the requirements that any accompanying (almost) arbitrary function must satisfy? The purpose of this work is to present the results of our investigations into these matters.\newline

We begin with the derivation of a simple form of the Master equation

\begin{equation} 
\int_{-\infty}^{\infty}\frac{F(a\,x(x+i))}{\cosh(\pi x)}dx=F(a/4) \label {master_basic},
\end{equation}

so-called because of its resemblance to what Hardy called Ramanujan's Master Theorem \cite{Berndt}
\begin{equation}
\int_0^{\infty}x^{s-1}F(x)dx=\Gamma(s)\phi(s)\,\,,
\label{RamMaster}
\end{equation}

where `` $\phi(s)$ is the `analytic continuation' of the Taylor coefficient of an unspecified function $F$''. Ramanujan \cite{Tichmarsh} applied \eqref{RamMaster} extensively on an {\it ad hoc} basis; subsequently, Hardy \cite{Hardy} identified sufficient conditions on $F$ for which it is valid. With the change of variables $u=cosh(\pi\,x)$, the identity \eqref{master_basic} resembles, but does not reduce to, the Berndt-Evans-Ramanujan \cite{Berndt} formula

\begin{equation}
\int_{-\infty}^{\infty}\frac{F(iw/u)}{1+u^2}du=F(w).
\label{BerEvRam}
\end{equation}

In addition to \eqref{RamMaster} and \eqref{BerEvRam}, the well-known Cauchy Integral formula forms a basis for other general ``Master'' equations, most of which involve rational functions \cite [Eq.021(12)] {GRHO}, a notable exception being \cite [Eq.335(17)] {GRHO}, which, with $a>0, b>0, F((z+1)/2)$ regular in $|z|\leq 1$, reads
\begin{equation}
\mapleinline{inert}{2d}{Int(f(exp(i*x)*cos(x))/(a^2*cos(x)^2+b^2*sin(x)^2),x = -Pi .. Pi) =
2*Pi*f(b/(a+b))/a/b;}{%
\[
{\displaystyle \int _{ - \pi }^{\pi }} {\displaystyle \frac {
\mathit{F}(e^{i\,x}\,\mathrm{cos}(x))}{a^{2}\,\mathrm{cos}(x)^{
2} + b^{2}\,\mathrm{sin}(x)^{2}}} \,dx={\displaystyle \frac {2\,
\pi \,}{a\,b}} \,\mathit{F}({\displaystyle  {b}/({a + b})} )\, .
\]
}
\label{Grob1}
\end{equation} 
In Section 2 we derive \eqref{master_basic} and a number of generalizations, and variations, thereby illustrating the technique used to discover them. We then proceed to quote related integration formulae obtained by applying these principles. In Section 3 and Appendix A, we exemplify the use of the results so-discovered by presenting some apparently new and interesting integration results. In Section 4, we conclude with a few observations and comments.\newline

\noindent
\section{Calculation} 

\noindent
 
Throughout, we shall use $z,x \in \mathbbm{C}$ where $z=a\,x\,(x+i)$ and assume that $F(z)$ is analytic  for $x\in\mathfrak{S}$ (defined below) and diverges less quickly than $\exp(\pi\,|x|)/x$ as $\mathfrak{Re}(x) \rightarrow \pm \infty$. In general, we let $a\in \mathbbm{R}$, although the possibility that $a\in \mathbbm{C}$ is sometimes considered in the discussion. In many cases however, we set $a=1$ without loss of generality. As well, we shall use $\{ j,k,n \} \in \mathbbm{N}$ except in some cases where these symbols $\in \mathbbm{Z}$ are allowed. The latter possibility is always specified where it applies. Define the complex strip 
\begin{equation}
\nonumber
\mathfrak{S} \equiv -1 \leq \mathfrak{Im}(x) \leq 0 \,.
\end{equation} 

{\bf Theorem}\vskip .1in

With the aforementioned provisos and notation, we have
\begin{equation}
\int_{-\infty}^{\infty}\frac{F(a\,x(x+i))}{\cosh(\pi x)}dx=F(a/4).
\label{Fcosh}
\end{equation}
\noindent
{\bf Proof}. Because $F(a\,x(x+i))$ is presumed analytic for $x\in\mathfrak{S}$, the only singularity of the integrand in this strip occurs when $x=-i/2$ (where $cosh(\pi x)$ vanishes) with residue $-(i/\pi)F(a/4)$. Therefore, after closing the contour at $\mathfrak{Re}(x) \rightarrow \pm \infty $, under the assumption that the integrand vanishes in these limits, we have 
\begin{equation}
\int_{-\infty}^{\infty}\frac{F(a\,x(x+i))}{\cosh(\pi x)}dx-\int_{-\infty-i}^{\infty-i}\frac{F(a\,x(x+i))}{\cosh(\pi x)}dx=2F(a/4). 
\label{Fdiff}
\end{equation}
But, with $x\rightarrow- x-i$ the second integral is the negative of the first. QED\vskip .1in

 By  a simple change of the variable of integration, one has\vskip .1in
 
 \noindent
 {\bf Corollary}\vskip .1in
 If $F(z)$ is analytic in the strip $-c \le \mathfrak{Re}(z)\le c$, where $0 \leq c <1/2$, then
\begin{equation}
\int_{c-i\infty}^{c+i\infty}\frac{F(a\,t(1-t)]}{\cos(\pi \,t)}\frac{dt}{2\pi\, i}=\frac{F(a/4)}{2\pi}. 
\label{c_Inf}
\end{equation}

\label{1.4}

\noindent
{\bf First Generalization}\vskip .1in

The derivation of \eqref{Fcosh} depends on the invariance (to within a multiplicative constant) of $\cosh(\pi\,x)$ under the transformation $x \rightarrow -x-i$ , and it is fair to ask if there exist any variants of the hyperbolic trigonometric functions that possess a similar property.
\newline

Consider a variant $K(x)$ as a candidate function to act as kernel, where
\begin{equation}
\mapleinline{inert}{2d}{K(x) = 1-(alpha*sinh(Pi*p*x)+beta*cosh(Pi*p*x))^2}{\[\displaystyle K \left( x \right) =1- \left( \alpha\,\sinh \left( \pi \,p\,x \right) +\beta\,\cosh \left( \pi \,p\,x \right)  \right) ^{2}\]} ,
\label{K(x)}
\end{equation}

temporarily requiring parameters $\alpha, \beta $ and $p \in \mathbbm{R}$.
\newline 

Apply the transformation $x \rightarrow -x-i$ to \eqref{K(x)} and rearrange the expansion of the result such that it reproduces the form \eqref{K(x)} with generalized (and complicated) coefficients. To reproduce the desired invariance, we require 

\begin{equation}
K(x)+K(-x-i)=0
\label{KplusK}
\end{equation}

\noindent
for all $x$, which in turn requires that 

\begin{equation}
\mapleinline{inert}{2d}{-2*alpha^2*cos(Pi*p)^2-2*beta^2*cos(Pi*p)^2+(4*I)*alpha*sin(Pi*p)*beta*cos(Pi*p) = 0}{\[\displaystyle -2\,{\alpha}^{2} \cos^{2} \left( \pi \,p \right)  -2\,{\beta}^{2} \cos ^{2}  \left( \pi \,p \right)   +4\,i\alpha\,\sin \left( \pi \,p \right) \beta\,\cos \left( \pi \,p \right) 
\mbox{}=0\]\, .}
\label{alphaEq}
\end{equation}

\noindent
The solutions are

\begin{equation}
\mapleinline{inert}{2d}{alpha = (I*tan(Pi*p)+I*sqrt(tan(Pi*p)^2+1))*beta}{\[\displaystyle \alpha= \left( \tan \left( \pi \,p \right) \pm\, \sqrt{ \left( \tan \left( \pi \,p \right)  \right) ^{2}+1} \right)\,i\, \beta\]} \,.
\label{alphaSoln}
\end{equation}

\noindent
Substituting \eqref{alphaSoln} into the expansion of \eqref{KplusK} gives a condition on $\beta$, depending on the sign of $\cos(\pi\,p)$

\begin{equation}
\mapleinline{inert}{2d}{-2*beta^2-2*beta^2*sin(Pi*p)+2*cos(Pi*p)^2 = 0}{\[\displaystyle -2\,{\beta}^{2}\,\pm\,2\,{\beta}^{2}\sin \left( \pi \,p \right) +2\,\cos ^{2} \left( \pi \,p \right) =0  .\]}
\label{betaEqn}
\end{equation} 

\noindent
The solutions of \eqref{betaEqn} are 

\begin{equation}
\mapleinline{inert}{2d}{cos(Pi*p)/sqrt(sin(Pi*p)+1)}{\[\displaystyle \beta = \pm {\frac {\cos \left( \pi \,p \right) }{ \sqrt{1\pm\,\sin \left( \pi \,p \right) }}}\]}
\label{betaSoln}
\end{equation}

\noindent
with a choice of sign depending on the sign of $\cos(\pi\,p)$ and the choice of solutions \eqref{alphaSoln}. After considering all possible sign permutations together with some algebraic simplification, the identification of \eqref{K(x)} subject to \eqref{KplusK} along with \eqref{alphaSoln} and \eqref{betaSoln} is

\begin{equation}
\mapleinline{inert}{2d}{K(x) = i*sinh(Pi*p*(2*x+i))}{\[\displaystyle K \left( x \right) = \pm \,i\sinh \left( \pi \,p \left( 2\,x+i \right)  \right) \]},
\label{G1}
\end{equation}

yielding the first generalization of \eqref{Fcosh}, subject to the same conditions that apply to \eqref{Fcosh} with the additional requirement that $p\in\mathbbm{R}$. Without loss of generality, let $p>0$. Then if $0<p<1$, the residue  at 
$x=-i/2$ of the result equivalent to \eqref{Fcosh} involving \eqref{G1} rather than $\cosh(\pi x)$ is easily determined, giving  

\begin{equation}
\mapleinline{inert}{2d}{Int(F(a*x*(x+i))/sinh(Pi*p*(2*x+i)), x = -infinity .. infinity) = -(1/2)*i*F((1/4)*a)/p}{\[\displaystyle \int _{-\infty }^{\infty }\!{\frac {F \left( ax \left( x+i \right)  \right) }{\sinh \left( \pi \,p \left( 2\,x+i \right)  \right) }}{dx}=-{\frac {i\,F \left( a/4 \right) 
\mbox{}}{2\,p}}\]}  
\label{gen1}
\end{equation}

reducing to \eqref{Fcosh} when $p=1/2$.
\newline 

In the case that $p>1$ it is easily shown that zeros of the reciprocal $sinh$ kernel in \eqref{gen1}, located at

\begin{equation}
\mapleinline{inert}{2d}{x_{p} = (1/2)*i*(1-p+n)/p}{\[\displaystyle x_{p}={\frac {i \left( 1-p+n \right) }{2\,p}}\]}
\label{x_p}
\end{equation}

where $n=\pm1,\pm2,\dots$ wander into the bounded integration strip $\mathfrak{S}$, and contribute their residues to the integral, if $-p<n\leq p-1$. In that case, (where $\lfloor p \rfloor \equiv $ the greatest integer less than $p$), the generalization of \eqref{gen1} becomes

\begin{equation}
\mapleinline{inert}{2d}{Int(F(a*x*(x+i))/sinh(2*Pi*p*x+I*p*Pi), x = -infinity .. infinity) = (1/2)*i*(Sum((-1)^n*F(-(1/4)*a*(n+1+p)*(n+1-p)/p^2), n = n1 .. n2))/p}{\[\displaystyle \int _{-\infty }^{\infty }\!{\frac {F \left( ax \left( x+i \right)  \right) }{\sinh \left( \pi \,p\,(2\,x+i ) \right) }}{dx}=i/(2p)
\mbox{}\,\sum _{n={\it \lfloor -p \rfloor}}^{{\it  \lfloor p-1 \rfloor}} \left( -1 \right) ^{n}F \left( -{\frac {a \left( (n+1)^2-p^2 \right) }{{4\,p}^{2}}} \right) \]\,.}
\label{gen1a}
\end{equation}

In the following Section a further generalization is obtained by removing the condition $p\in\mathbbm{R}$. 
\newline

{\bf Second Generalization} \vskip 0.1 in

In \eqref{K(x)} let
 
\begin{equation}
p=q+i\,r
\label{p=q+i*r}
\end{equation}

with $q,r \in \mathbbm{R}$. Direct substitution of \eqref{p=q+i*r} into \eqref{G1} shows that \eqref{KplusK} remains satisfied under this generalization and it only remains to locate the travelling zeros of \eqref{G1} for a range of $q$ and $r$. In exact analogy to \eqref{x_p}, the zeros of \eqref{G1} are located at the points

\begin{equation}
\mapleinline{inert}{2d}{x[p] = (1/2)*n*r/(r^2+q^2)+(1/2)*i*(-1+n*q/(r^2+q^2))}{\[\displaystyle x_{{p}}=1/2\,\left( {\frac {nr}{{r}^{2}+{q}^{2}}}\right) +i/2 \left( -1+{\frac {nq}{{r}^{2}+{q}^{2}}} \right) \]}
\label{xp}
\end{equation}

where, again, $n=\pm1,\pm2,\dots$ . In order that any zero contributes its corresponding residue to the integral, we require that 

\begin{equation}
\mapleinline{inert}{2d}{n < (r^2+q^2)/q, -(r^2+q^2)/q < n}{\[\displaystyle -{\frac {{r}^{2}+{q}^{2}}{q}} < n<{\frac {{r}^{2}+{q}^{2}}{q}} \]}
\label{nlims}
\end{equation}

in the case $q>0$, with the inequality signs reversed if $q<0$. The generalization of \eqref{gen1a} becomes

\begin{equation}
\mapleinline{inert}{2d}{Int(F(a*x*(x+i))/sinh((q+i*r)*Pi*(i+2*x)),x = -infinity .. infinity)
=
-1/2*i*F(1/4*a)/(q+i*r)-i*sum((-1)^n*F(1/4*a*(n+q+i*r)*(-n+q+i*r)/(q+i
*r)^2),n = 1 .. M1)/(q+i*r);}{%
\maplemultiline{
{\displaystyle \int _{ - \infty }^{\infty }} {\displaystyle 
\frac {\mathit{F}(a\,x\,(x + i))}{\mathrm{sinh}(\pi 
\,(q + i\,r)\,(i + 2\,x))}} \,dx= 
 - {\displaystyle \frac {i}{2}} \,{\displaystyle \frac {\,
\mathit{F}({\displaystyle a/4} )}{q + i\,r}}  - 
{\displaystyle \frac{i}{q + i\,r}  {\,\,\! {\displaystyle \sum _{n=1}^{
\mathit{M}}} \,(-1)^{n}\,\mathit{F}\left({\displaystyle \frac {a\,(  (q + i\,r)^2 - n^2)}{4\,(q + i\,r)^{2}}} \right ) \! 
 }}  }
}
\label{gen2}
\end{equation}

where $M=\lfloor \pm\frac{q^2+r^2}{q} \rfloor$, according as $q>0$ or $q<0$. In the case $q=0$, carefully following the same method eventually yields the limiting case

\begin{equation}
\mapleinline{inert}{2d}{i*r*Int(F(a*x*(x+i))/sin(r*Pi*(i+2*x)),x = -infinity .. infinity) =
Sum((-1)^n*F(1/4*a*(r^2+n^2)/r^2),n = 1 .. infinity)+1/2*F(1/4*a);}{%
\[
i\,r\,{\displaystyle \int _{ - \infty }^{\infty }} 
{\displaystyle \frac {\mathit{F}(a\,x\,(x + i))}{\mathrm{sin}(r\,
\pi \,(i + 2\,x))}} \,dx=  \! {\displaystyle \sum _{n=1}^{
\infty }} \,(-1)^{n}\,\mathit{F}\left({\displaystyle \frac {a\,(r^{2}
 + n^{2})}{4\,r^{2}}} \right) \!   + {\displaystyle \frac {1}{2
}} \,\mathit{F}({\displaystyle a/4} )\,\,.
\]
}
\label{q=0}
\end{equation}
\newline

Further generalization is possible by noting that odd powers of \eqref{G1} obey \eqref{KplusK},
giving for example

\begin{equation}
\mapleinline{inert}{2d}{Case3c := Int(F(a*x*(x+i))/sinh((q+i*r)*(i+2*x))^3,x = -infinity ..
infinity) = 1/4*i*Pi*F(1/4*a)/(q+i*r)-1/8*i*a*eval(diff(F(x),x),\{x =
1/4*a\})/Pi/(q+i*r)^3+1/2*i*Pi*Sum((-1)^n*F(1/4*A[n]),n = 1 ..
M1)/(q+i*r)-1/4*i*a*Sum((-1)^n*eval(diff(F(x),x),\{x = 1/4*A[n]\}),n =
1 ..
M1)/Pi/(q+i*r)^3+1/8*i*a^2*Sum((-1)^n*eval(diff(F(x),`$`(x,2)),\{x =
1/4*A[n]\})*n^2,n = 1 .. M1)/Pi/(q+i*r)^5;}{%
\maplemultiline{
{\displaystyle \pi\,\int _{ - \infty }^{\infty }} 
{\displaystyle \frac {\mathit{F}(a\,x\,(x + i))}{\mathrm{sinh}^{3} (\pi\,(q
 + i\,r)\,(i + 2\,x)) }} \,dx={\displaystyle \frac {i\,\pi \,}{4}} \,
{\displaystyle \frac {\mathit{F}({\displaystyle {a/4}} )}{(q + i\,r)}}  - {\displaystyle \frac {i\,a\,}{8\,\pi \,(q + i\,r)^{3}}} \,
{\displaystyle {{ {\frac {d}{dx}}\,\mathit{F}(x)
\vrule \lower 11pt \hbox{$\displaystyle \, x={\displaystyle 
{a/4}} $}}}}  \\
\mbox{} + {\displaystyle \frac {i
\,\pi}{2\,(q + i\,r)}} \,{\displaystyle {  \! {\displaystyle \sum _{n=1}^{\mathit{M}}} \,(
-1)^{n}\,\mathit{F}({\displaystyle {x_{n}\,/4}} \,{}) \! 
 }}  - {\displaystyle \frac {i\,a}{4\,\pi \,(q + i\,r)^{3}}} \,
{\displaystyle {  {\displaystyle \sum _{n=1
}^{\mathit{M}}} \,(-1)^{n}\,{ {\frac {d}{dx}}\,\mathit{F}(x)
\vrule \lower 5pt \hbox{$ \, x={x_{n}/4}$}} \! 
 }}  \\
\mbox{} + {\displaystyle \frac {i\,a^{2}}{8\,\pi \,(q + i\,r)^{5}}} \,{\displaystyle  {
\ {\displaystyle \sum _{n=1}^{\mathit{M}}} \,
(-1)^{n}\,n^{2} { {\frac {d^{2}}{dx^{2}}}\,\mathit{F}(x)
\vrule \lower 5pt \hbox{$ \, x={x_{n}/4}$}}\, \! }{}}  }
\label{cubed}
}
\end{equation}

where M is defined in \eqref{gen2} and
\begin{equation} 
\mapleinline{inert}{2d}{A[n] = a*(-n+q+i*r)*(n+q+i*r)/(q+i*r)^2}{\[\displaystyle x_{{n}}={\frac {a \left( -n+q+ir \right)  \left( n+q+ir \right) }{ \left( q+ir \right) ^{2}}}\]\,.}
\label{A[n]}
\end{equation} 

In contrast, even powers of \eqref{G1} obey an analogue of \eqref{KplusK}, namely 

\begin{equation}
K(x)-K(-x-i)=0\,.
\label{KminusK}
\end{equation}

This means that the right hand side of the equivalent of \eqref{Fdiff} vanishes. However, as \eqref{x_p} attests, the integrand still contains poles in $\mathfrak{S}$ for certain values of $p$. For the case $n=0$ the residue of the pole at $x=-i/2$ vanishes. For other values of $n$, the poles occur in pairs symmetrically distributed about the horizontal line passing through $x=-i/2$, and their residues have opposite sign. In that case, the sum of the residues vanish pairwise, consistent with the right-hand side of \eqref{Fdiff} vanishing.\newline 

The fact that integrals possessing kernels which simultaneously involve even powers of the denominator and satisfy \eqref{KplusK} evade discovery, suggests that the lower boundary of the integration region be shifted to pass through the point $x=-i/2$ rather than $x=-i$ to destroy the (inconvenient) symmetry. When this is done, a special case closely related to the results already discussed here is discovered:

\begin{equation}
K_1(x) = {\displaystyle \sinh(2\,k\,\pi\,(i+2\,x))=\sinh(4\,k\,\pi \,x)}
\end{equation}

which obeys
\begin{equation}
K_{1}(x)+K_{1}(-x-i/2)=0
\label{K with i/2}.
\end{equation}
where $k=0,\pm1,\pm2,\dots. $
\newline

Following the methods developed above, and, with a minor change of integration variables we find the corresponding result

\begin{equation}
\mapleinline{inert}{2d}{Int(F(a*x*(i+x))/sinh(2*k*Pi*(i+x)),x = -infinity ..
infinity)+Int(F(a*x*(i+x))/sinh(2*k*Pi*x),x = -infinity .. infinity) =
-i*(-1)^k*Sum(F(1/4*a*(k-n)*(n+k)/k^2)*(-1)^n,n = 1 ..
k-1)/k-1/2*i*(-1)^k*F(1/4*a)/k;}{%
\maplemultiline{
{\displaystyle \int _{ - \infty }^{\infty }} {\displaystyle 
\frac {\mathit{F}(a\,x\,(i + x))}{\mathrm{sinh}(2\,k\,\pi \,(i + 
x))}} \,dx = {\displaystyle \int _{ - \infty }^{\infty }} 
{\displaystyle \frac {\mathit{F}(a\,x\,(i + x))}{\mathrm{sinh}(2
\,k\,\pi \,x)}} \,dx= 
 - {\displaystyle \frac {i\,(-1)^{k}}{k}\,} \left ( {\displaystyle  {\displaystyle 
\sum _{n=1}^{k - 1}} \,(-1)^{n}\mathit{F}({\displaystyle{x_{n}}} )\, \!   }  + 
{\displaystyle \frac {1}{2}} \,{\displaystyle 
\,\mathit{F}({\displaystyle \frac {a}{4}} )   }\right ) }
}
\label{A1}
\end{equation}

where
\begin{equation}
x_{n}=\frac{a\,(k^2 - n^2)}{4\,k^2}
\end{equation}

and $a\in \mathbbm{C}$ must lie in a range where $F(z)\in \mathfrak{S}$ and must vanish faster than $O(1)$ as $x \rightarrow 0$. Setting $p=k$ in \eqref{gen1a} corresponds to \eqref{A1} up to a factor $(-1)^k$ on the right-hand side.
\newline

In a similar vein, another related result may be obtained by bounding the lower integration boundary at $-x-(2\,k-1)\,i)$ rather than $-x-i$ and after some algebra the following generalization, which reduces to \eqref{Fcosh} when $k=1$, and corresponds (up to a factor of $i (-1)^{k+1}$ ) to \eqref{gen1a} when $p=k-1/2$ is found 

\begin{equation}
\mapleinline{inert}{2d}{Int(F(a*v*(v+i))/cosh(Pi*(2*k-1)*v),v = -infinity .. infinity) =
-2*(-1)^k/(2*k-1)*Sum((-1)^n*F(1/4*a*(2*k-2*n-1)*(2*k+2*n-1)/(2*k-1)^2
),n = 1 .. k-1)-(-1)^k*F(1/4*a)/(2*k-1);}{%
\maplemultiline{
{\displaystyle \int _{ - \infty }^{\infty }} {\displaystyle 
\frac {\mathit{F}(a\,x\,(x + i))}{\mathrm{cosh}(\pi \,(2\,k - 1)
\,x)}} \,dx=  {\displaystyle {(-1)^{(k+1)}\,
\frac {\mathit{F}({\displaystyle a/4} )}{(2\,k - 1)}}{}} 
   -{\displaystyle {\frac {2\,(-1)^{k}}{2\,k - 1}\,\,\, {\displaystyle 
\sum _{n=1}^{k - 1}} \,(-1)^{n}\,\mathit{F}\left({\displaystyle 
\frac {a\,((2\,k-1)^2 - 4\,n^2)}{4\,(2\,k - 1)^{2}
}} \right) }{}} \,.  }
}
\label{A.3}
\end{equation}

All the above have been written in a canonically similar fashion by redefining variables of integration, pairing residues where possible, and utilizing the convention that any sum vanishes when a lower limit exceeds an upper limit.
\newline

Finally, another result based on the same principles, but lacking some of the symmetry of the foregoing has been identified:
\begin{equation}
\mapleinline{inert}{2d}{A4a :=
Int(F(a*x*(x+i))*cosh(Pi*(2*k-1)*x)/(cosh(Pi*b)^2+sinh(Pi*(2*k-1)*x)^2
),x = -infinity .. infinity) =
-Sum((-1)^n*F(-1/4*i*a*(2*i*b-2*n-1+4*k)*(i*(2*n-1)+2*b)/(2*k-1)^2),n
= 1 .. 2*k-1)/cosh(Pi*b)/(2*k-1);}{%
\maplemultiline{
{\displaystyle \int _{ - \infty }^{\infty }} 
{\displaystyle \frac {\mathit{F}(a\,x\,(x + i))\,\mathrm{cosh}(
\pi \,(2\,k - 1)\,x)}{\mathrm{cosh^{2}}(\pi \,b) + \mathrm{sinh^{2}}(
\pi \,(2\,k - 1)\,x)}} \,dx =\mapleinline{inert}{2d}{Int(F(a*x*(i+x))*cosh(Pi*(2*k-1)*x)/cosh(Pi*(b+(2*k-1)*x))/cosh(Pi*(b
-(2*k-1)*x)),x = -infinity .. infinity);}{%
\[
{\displaystyle \int _{ - \infty }^{\infty }} {\displaystyle 
\frac {\mathit{F}(a\,x\,(i + x))\,\mathrm{cosh}(\pi \,(2\,k - 1)
\,x)}{\mathrm{cosh}(\pi \,(b + (2\,k - 1)\,x))\,\mathrm{cosh}(\pi
 \,(b - (2\,k - 1)\,x))}} \,dx
\]
}
\\  
= {\displaystyle \frac{1}{\mathrm{cosh}(\pi \,b)\,(2\,k - 1)}\displaystyle  {{\displaystyle \sum _{n=1}^{2\,k - 1}} 
\,(-1)^{n+1}\,\mathit{F}\left( - {\displaystyle \frac {i\,a\,(2\,i\,b - 
2\,n - 1 + 4\,k)\,(i\,(2\,n - 1) + 2\,b)}{4\,(2\,k - 1)^{2}}} \right)}{
}} \, . }
\label{A.4}
}
\end{equation}

If $b=0$, \eqref{A.4} reduces to \eqref{A.3}. Note that (\ref{A.4}) is valid for $b \in \mathbbm{C}$ where the integral exists, and that the right-hand side provides an analytic continuation for those regions of $b$ where the integral does not exist.
\newline

{\bf Variations} \vskip 0.1 in

In any of the above, the function $F(ax(x+i))$ can be replaced by products or sums of specialized functions, provided that such product or sum is invariant under $x \rightarrow -x-i$. For example, if a function $F_{o|e}(x)$ is odd or even under $x\rightarrow -x$, the product $F_{o|e}(ax) F_{o|e}(a(x+i))$ (but not  $F_{o|e}(ax) F_{e|o}(a(x+i))$) can replace $F(ax(x+i))$ in the integrand of any of the above and any appearance of the form $F(z_{1} z_{2})$ in the right-hand side of any of the results quoted previously, becomes $F_{o|e}(z_{1}) F_{o|e}(z_{2})$. This leads, for example, to the following variation of \eqref{gen2}

\begin{equation}
\mapleinline{inert}{2d}{Int(F_[oe](a*x)*F_[oe](a*(x+i))/sinh((q+i*r)*Pi*(i+2*x)),x =
-infinity .. infinity) =
-1/2*i*F_[oe](-1/2*i*a)/(q+i*r)*F_[oe](1/2*i*a)-1/2*i*sum(2*F_[oe](-1/
2*i*(q+i*r-n)*a/(q+i*r))*F_[oe](1/2*i*(q+i*r+n)*a/(q+i*r))*(-1)^(-n),n
= 1 .. M1)/(q+i*r);}{%
\maplemultiline{
{\displaystyle \int _{ - \infty }^{\infty }} {\displaystyle 
\frac {{\mathit{F} _{ {o | {e} }}}(a\,x)\,{\mathit{F}_{
\mathit{o|e}}}(a\,(x + i))}{\mathrm{sinh}(\pi \,(q + i\,r)\,(i + 2
\,x))}} \,dx= -\,{\displaystyle \frac {i}{2}} \,{\displaystyle 
\frac {{\mathit{F}_{\mathit{o|e}}}( - {\displaystyle \frac {i
\,a}{2}} )\,{\mathit{F}_{\mathit{o|e}}}({\displaystyle \frac {i
\,a}{2}} )}{q + i\,r}}  \\
\mbox{} + {\displaystyle \frac{i}{(q + i\,r)}
\, \! \displaystyle} \,{\displaystyle { {\displaystyle \sum _{n=1}^{\mathit{M}}} \,{
\,(
-1)^{ n+1}\mathit{F}_{\mathit{o|e}}}( - {\displaystyle \frac {i\,(q + i\,r
 - n)\,a}{2\,(q + i\,r)}} )\,{\mathit{F}_{\mathit{o|e}}}(
{\displaystyle \frac {i\,(q + i\,r + n)\,a}{2\,(q + i\,r)}} ) \!  }{}}  }
}
\label{var1}
\end{equation}

where $M$ is the same as before. Similarly, in any of the above, we may alternatively replace $F(ax(x+i))$ by either a sum of even functions $F_{e}(ax)+F_{e}(a(x+i))$ or a difference of odd functions $F_{o}(ax)-F_{o}(a(x+i))$, motivated by the knowledge that any function $F(x)$ can always be rewritten as a sum of $F_{e}(x)$ and $F_{o}(x)$. To clarify how the corresponding parts of the right hand side are written, a number of variations of \eqref{gen1a} are listed explicitly in Appendix A.

\section{Examples}

\subsection{Example: $F(z)=1$}

The simplest example comes about by setting $F(z)=1$ in \eqref{gen1a}. By utilizing the symmetry of the integrand, splitting it into real and imaginary parts and applying some elementary simplification, we find
\begin{equation}
\mapleinline{inert}{2d}{8*p*sin(Pi*p)*Int(cosh(2*Pi*p*x)/(cos(2*Pi*p)-cosh(4*Pi*p*x)),x = 0
.. infinity) = -1/2*(-1)^floor(p)+1/2*(-1)^floor(-p);}{%
\[
4\,p\,\,\mathrm{sin}(\pi \,p)\,{\displaystyle \int _{0}^{\infty }} 
{\displaystyle \frac {\mathrm{cosh}(\pi \,p\,t)}{\mathrm{cos}(
2\,\pi \,p) - \mathrm{cosh}(2\,\pi \,p\,t)}} \,dt=  
\frac{1}{2}\left( (-1)^{\lfloor -p \rfloor}-(-1)^{\lfloor p \rfloor})\right) 
\]
}
\label{X1}
\end{equation}

valid for all $p>0$. This example has the interesting property that the integrand contains only continuous, and non-pathological functions with a simple non-integrable singularity at the origin (for distinct values of the parameter $p$) which can be treated by interpreting the integral in the principal value sense, yet the left-hand side itself represents an analytically pathological square wave of magnitude unity and period two as a function of $p$, provided only that $p \neq k$. In \eqref{X1}, both sides vanish when $p=k, k\neq 0$. With reference to \eqref{A1}, this example violates the requirement that $F(z)\rightarrow 0$ as $x \rightarrow 0$; since the right-hand side of that result offers $1/2$ as a representation of \eqref{X1} when $p=k$, this suggests that \eqref{A1} needs to be recast as a principal value integral in order to apply this example to that result. It is worth noting that both Maple \cite{Maple} and Mathematica \cite{Math} arrive at an equivalent, but less elegant, result for this case. 
\newline
\subsection{Example: $F(z)=z^{j}$}
A second example is obtained by setting $F(z)=z^j$ with $j=1,2,\dots$ in \eqref{gen1a}. Comparing the real and imaginary parts and applying symmetry gives two general results for even and odd values of $j$, neither of which appears to be known to either Maple \cite{Maple} or Mathematica \cite{Math} for general values of $j$ and $p$, although both codes can evaluate special cases. In the first case we set $j \rightarrow 2j$; in the second we have set $j \rightarrow 2j-1=k$ to arrive at the following: 

\begin{equation}
\mapleinline{inert}{2d}{X2a :=
Int(t^(2*j)*(t^2+1)^j*(sin(2*j*arctan(t))*sinh(2*Pi*p*t)*cos(Pi*p)+cos
(2*j*arctan(t))*cosh(2*Pi*p*t)*sin(Pi*p))/(cos(Pi*p)^2-cosh(2*Pi*p*t)^
2),t = 0 .. infinity) =
1/4*(-1)^j*Sum((-1)^n*(1/4*(n+1+p)*(n+1-p)/p^2)^{2*j},n = floor(-p)
.. -1+floor(p))/p;}{%
\maplemultiline{
{\displaystyle \int _{0}^{\infty }} 
{\displaystyle \frac {t^{2j}\,(t^{2}+ 1)^{j}\,(\mathrm{sin}(
2\,j\,\mathrm{arctan}(t))\,\mathrm{sinh}(2\,\pi \,p\,t)\,\mathrm{
cos}(\pi \,p) + \mathrm{cos}(2\,j\,\mathrm{arctan}(t))\,\mathrm{
cosh}(2\,\pi \,p\,t)\,\mathrm{sin}(\pi \,p))}{\mathrm{cos^{2}}(\pi \,
p) - \mathrm{cosh^{2}}(2\,\pi \,p\,t)}} dt  \\
={\displaystyle } \,{\displaystyle \frac {(-1)^{j}} {2\,(2\,p)^{4\,j+1}}
 \! {\displaystyle \sum _{n= \lfloor - p\rfloor }^{ - 1
 + \lfloor p \rfloor }} \,(-1)^{n}\,{\displaystyle {((n + 1)^2 - p^2)}{}}^{ 2\,j} } 
 }
}
\label{x2a}
\end{equation}

\begin{equation}
\mapleinline{inert}{2d}{Int(t^k*(t^2+1)^(1/2*k)*(cos(k*arctan(t))*sinh(2*Pi*p*t)*cos(Pi*p)-si
n(k*arctan(t))*cosh(2*Pi*p*t)*sin(Pi*p))/(cos(Pi*p)-cosh(2*Pi*p*t))/(c
os(Pi*p)+cosh(2*Pi*p*t)),t = 0 .. infinity) =
1/4*(-1)^j*Sum((-1)^n*(-1/4*(n+1+p)*(n+1-p)/p^2)^k,n = floor(-p) ..
-1+floor(p))/p;}{%
\maplemultiline{
{\displaystyle \int _{0}^{\infty }} { \displaystyle  \frac { t^{k}\,(t^{2} + 1)^{k/2} \\
(\mathrm{cos}(k\,\mathrm{arctan}(t))\,\mathrm{sinh}(2\,\pi \,p\,t
)\,\mathrm{cos}(\pi \,p) - \mathrm{sin}(k\,\mathrm{arctan}(t))\,
\mathrm{cosh}(2\,\pi \,p\,t)\,\mathrm{sin}(\pi \,p))} { \displaystyle \mathrm{cos^{2}}(\pi \,
p) - \mathrm{cosh^{2}}(2\,\pi \,p\,t)}} dt \\
={\displaystyle \frac{ {(-1)^{k+\lfloor k/2 \rfloor } }}
{ 2\,(2\,p)^{2\,k+1}}  \! {\displaystyle \sum _{n= \lfloor - p \rfloor}^{ - 1
 + \lfloor p \rfloor}} \,(-1)^{n}\,{\displaystyle  {((n + 1)^2 - p^2)}{}} ^{k} \!   }{}  }
}
\label{x2b}
\end{equation}

The obviously interesting special cases correspond to $p=k-1/2$ and $p=k$. In the first and second instances respectively, from \eqref{x2a} we find

\begin{equation}
\mapleinline{inert}{2d}{X2a1 := Int(t^(2*j)*(t^2+1)^j*cos(2*j*arctan(t))/cosh(Pi*(2*k-1)*t),t
= 0 .. infinity) =
1/2*(-1)^(j+k)*sum((-1)^n*((2*n+1+2*k)*(-2*n-3+2*k))^(2*j),n =
floor(-k+1/2) .. -2+floor(k+1/2))/(16^j)/((2*k-1)^(4*j+1));}{%
\maplemultiline{
{\displaystyle \int _{0}^{\infty }} 
{\displaystyle \frac {t^{2\,j}\,(t^{2} + 1)^{j}\,\mathrm{cos}(2
\,j\,\mathrm{arctan}(t))}{\mathrm{cosh}(\pi \,(2\,k - 1)\,t)}} \,
dt  =
 \,{\displaystyle \frac {(-1)^{(j + k
)}\,{\displaystyle \sum _{n= - k }^{ k - 2}} \,(-1)^{n}\,((2\,n + 1 + 2\,k)\,( - 2\,
n - 3 + 2\,k))^{2\,j}}{(4\,k - 2)^{(4\,j + 1)}}}  }
}
\label{x2a1}
\end{equation}
and
\begin{equation}
\mapleinline{inert}{2d}{X2a2 := Int(t^(2*j)*(t^2+1)^j*sin(2*j*arctan(t))/sinh(2*Pi*k*t),t = 0
.. infinity) =
1/4*(-1)^(j+k+1)*(1/16)^j/(k^(4*j+1))*sum((-1)^n*((n+1+k)*(-n-1+k))^(2
*j),n = floor(-k) .. -1+floor(k));}{%
\maplemultiline{
{\displaystyle \int _{0}^{\infty }} 
{\displaystyle \frac {t^{2\,j}\,(t^{2} + 1)^{j}\,\mathrm{sin}(2
\,j\,\mathrm{arctan}(t))}{\mathrm{sinh}(2\,\pi \,k\,t)}} \,dt
={\displaystyle } \,{\displaystyle \frac {(-1)^{j + k
 + 1}\,\,{\displaystyle 
\sum _{n= - k }^{ k - 1  }} \,(-1)
^{n}\,(k^2 -(n + 1)^2)^{2\,j}}{2\,(2\,k)^{4\,j + 1}}.} 
 }
}
\label{x2a2}
\end{equation}

For the corresponding cases based on \eqref{x2b} we find

\begin{equation}
\mapleinline{inert}{2d}{X2b1 :=
Int(t^(2*j-1)*(t^2+1)^(j-1/2)*sin((2*j-1)*arctan(t))/cosh(Pi*(2*k-1)*t
),t = 0 .. infinity) =
-2*(-1)^(j+k)*sum((-1)^n*((2*n+1+2*k)*(-2*n-3+2*k))^(2*j)/(2*n+1+2*k)/
(-2*n-3+2*k),n = -k .. -2+k)/(16^j)/((2*k-1)^(4*j-1));}{%
\maplemultiline{
{\displaystyle \int _{0}^{\infty }} 
{\displaystyle \frac {t^{2\,j - 1}\,(t^{2} + 1)^{j - \frac{1}{2} }\,
\mathrm{sin}((2\,j - 1)\,\mathrm{arctan}(t))}{\mathrm{cosh}(\pi 
\,(2\,k - 1)\,t)}} \,dt  =
 {\displaystyle \frac {(-1)^{j + k + 1}\, 
{\displaystyle \sum _{n= - k}^{k - 2}} \,{\displaystyle 
 {(-1)^{n}\,((2\,n + 1 + 2\,k)\,( - 2\,n - 3 + 2\,k))^{2\,j
-1}}{}} }{\,(4\,k - 2)^{4\,j - 1}}}  }
}
\label{x2b1}
\end{equation}

\begin{equation}
\mapleinline{inert}{2d}{X2b2 :=
Int(t^(2*j-1)*(t^2+1)^(j-1/2)*cos((2*j-1)*arctan(t))/sinh(2*Pi*k*t),t
= 0 .. infinity) =
(-1)^(j+k+1)*16^(-j)*k^(1-4*j)*sum((-1)^n*((n+1+k)*(-n-1+k))^(2*j)/(n+
1+k)/(-n-1+k),n = -k .. -1+k);}{%
\maplemultiline{
{\displaystyle \int _{0}^{\infty }} 
{\displaystyle \frac {t^{2\,j - 1}\,(t^{2} + 1)^{j - 1/2}\,
\mathrm{cos}((2\,j - 1)\,\mathrm{arctan}(t))}{\mathrm{sinh}(2\,
\pi \,k\,t)}} \,dt  =
(-1)^{(j + k + 1)}\, 
\frac{{\displaystyle \sum _{n= - k}^{ k - 1 }} \,{\displaystyle 
{(-1)^{n}\,( k^2-( n + 1)^2)^{(2\,j-1)}}}  }{\displaystyle 2\,(2\,k)^{4\,j-1}}}.
}
\label{x2b2}
\end{equation}
\newline

It is worth noting that both \eqref{x2a} and \eqref{x2b} are intriguingly similar to integral representations of the Dirichlet eta function $\eta(2j)$ and $\eta(2j-1)$ respectively (see \eqref{eta} below), that the numerator trigonometric functions reduce to a rational function in $t^2$ (\cite{Luke}, page 210, Eqs. 6.2.1(11) and 6.2.1(12)) and that in the case $j=1$, \eqref{x2b} approximates a sawtooth wave as a function of $p$.
\newline
\subsection{Example: $F(z)=1/z^{s}$}
Carrying the previous example further, in \eqref{A.3} let
\begin{equation}
F(z)=1/z^s \,;
\end{equation}

furthermore allow $s \in \mathbbm{C}$ such that $\Re(s)<1$. Even though $F(z)$ is cut when $\Re (z)<0, \Re(s) \not \in \mathbb{Z} $, the logic of the derivation follows as before, and we find

\begin{equation}
\mapleinline{inert}{2d}{Rk :=
Int(t^(-s)*(t^2+1)^(-1/2*s)*cos(s*arctan(1/t))/cosh(Pi*(2*k-1)*t),t =
0 .. infinity) =
2^s*(-1)^(1+k)*(2*k-1)^(s-1)*Sum((-1)^n*(1/(2*k+2*n-1)-1/(-2*k+2*n+1))
^s,n = 1 .. k-1)-2^(2*s-1)*(-1)^k/(2*k-1);}{%
\maplemultiline{
{\displaystyle \int _{0}^{\infty }} 
{\displaystyle \frac {t^{-s}\,(t^{2} + 1)^{ -s/2 }
\,\mathrm{cos}(s\,\mathrm{arctan}({\displaystyle 1/t } ))
}{\mathrm{cosh}(\pi \,(2\,k - 1)\,t)}} \,dt \\
=2^{2\,s}\,(-1)^{k+1}\,(2\,k - 1)^{2\,s - 1}\, 
{\displaystyle \sum _{n=1}^{k - 1}} \,(-1)^{n}\,\left ( {\displaystyle 
\frac {1}{(2\,k-1)^2 - 4\,n^2}}  \right ) ^{s} \!  
\mbox{} + {\displaystyle \frac {2^{(2\,s - 1)}\,(-1)^{k+1}}{2\,k - 
1}} \,. }
}
\label{X2a}
\end{equation}
\newline

It is interesting to compare the left-hand side of \eqref{X2a} to one of the (known) integral representations \cite [Eq. (1.12(15))] {Bateman} of the Riemann zeta function $\zeta(s)$:

\begin{equation}
\mapleinline{inert}{2d}{eta(s) =
2^(s-1)*Int(cos(s*arctan(t))(t^2+1)^(-1/2*s)/cosh(1/2*Pi*t),t = 0 ..
infinity);}{%
\[
\zeta (s)={\displaystyle \frac {2^{s - 1}}{1-2^{1-s}}}\,{\displaystyle \int _{0}^{\infty }} 
{\displaystyle \frac {\mathrm{cos}(s\,\mathrm{arctan}(t))(t^{2}
 + 1)^{( - s/2 )}}{\mathrm{cosh}({\displaystyle \pi\,t/2  } )}} \,dt \, .
\]
}
\label{eta}
\end{equation}

Further integral evaluations can be obtained by either comparing real and imaginary parts of \eqref{X2a} for complex values of $s$, or by utilizing the well-known result

\begin{equation}
{\mathrm{arctan}} (v) +{\mathrm {arctan}}(1/v)=\pi/2 , \,\,\; v>0
\label{arctan}
\end{equation}

which yields the identities

\begin{equation}
\mapleinline{inert}{2d}{sin(2*s*arctan(1/v)) =
sin(s*Pi)*cos(2*s*arctan(v))-cos(s*Pi)*sin(2*s*arctan(v));}{%
\[
\mathrm{sin}(2\,s\,\mathrm{arctan}({\displaystyle 1/v} )
)=\mathrm{sin}(s\,\pi )\,\mathrm{cos}(2\,s\,\mathrm{arctan}(v))
 - \mathrm{cos}(s\,\pi )\,\mathrm{sin}(2\,s\,\mathrm{arctan}(v))
\]
}
\end{equation}
\begin{equation}
\mapleinline{inert}{2d}{cos(2*s*arctan(1/v)) =
cos(s*Pi)*cos(2*s*arctan(v))+sin(s*Pi)*sin(2*s*arctan(v));}{%
\[
\mathrm{cos}(2\,s\,\mathrm{arctan}({\displaystyle 1/v} )
)=\mathrm{cos}(s\,\pi )\,\mathrm{cos}(2\,s\,\mathrm{arctan}(v))
 + \mathrm{sin}(s\,\pi )\,\mathrm{sin}(2\,s\,\mathrm{arctan}(v))\,,
\]
\label{SinCosArctan}
}
\end{equation}

structures that persistently appear in the examples that follow.
\newline

\subsection{Example $F_{o|e}(az)=\cos(az)$}
With respect to the variation \eqref{var1} , letting $F_{o|e}(az)=cos(az)$, repeating the same steps as before (with $r=0$), and using Maple \cite{Maple} to explicitly sum the right-hand side, leads to

\begin{equation}
\mapleinline{inert}{2d}{sin(Pi*q)*cosh(a)*Int((cos(2*a*v)+1)*cosh(2*Pi*q*v)/(-cos(Pi*q)^2+cos
h(2*Pi*q*v)^2),v = 0 ..
infinity)+cos(Pi*q)*sinh(a)*Int(sin(2*a*v)*sinh(2*Pi*q*v)/(-cos(Pi*q)^
2+cosh(2*Pi*q*v)^2),v = 0 .. infinity) =
1/8*(2*cosh((a*floor(q)+a)/q)+2*cosh(a*floor(q)/q)+cosh((a*q-a)/q)+cos
h((a*q+a)/q)+2*cosh(a))/((-1)^floor(q))/q/(cosh(a/q)+1);}{%
\maplemultiline{
2\,\mathrm{sin}(\pi \,q)\,\mathrm{cosh}(a)\,{\displaystyle \int _{0}
^{\infty }} {\displaystyle \frac {\mathrm{cos^2}(a\,v)\,
\mathrm{cosh}(2\,\pi \,q\,v)}{\mathrm{cosh}^2(2\,\pi \,q\,v) - \mathrm{cos}^2(\pi \,q)  
}} \,dv \\
\mbox{} + \mathrm{cos}(\pi \,q)\,\mathrm{sinh}(a)\,
{\displaystyle \int _{0}^{\infty }} {\displaystyle \frac {
\mathrm{sin}(2\,a\,v)\,\mathrm{sinh}(2\,\pi \,q\,v)}{\mathrm{cosh}^2(2\,\pi \,q\,v) - \mathrm{
cos}^2(\pi \,q) }} \,dv\\
={\displaystyle } \,{\displaystyle \frac {2\,\mathrm{
cosh}({\displaystyle \frac {a\,\mathrm{\lfloor}q\rfloor + a}{q}} ) + 2\,
\mathrm{cosh}({\displaystyle \frac {a\,\mathrm{\lfloor}q \rfloor}{q}} )
 + \mathrm{cosh}({\displaystyle \frac {a\,q - a}{q}} ) + \mathrm{
cosh}({\displaystyle \frac {a\,q + a}{q}} ) + 2\,\mathrm{cosh}(a)
}{8\,q\,(-1)^{\mathrm{\lfloor q\rfloor}}\,(\mathrm{cosh}({\displaystyle 
\frac {a}{q}} ) + 1)}}  }
}
\label{case5a}
\end{equation}

which appears to be reducible to combinations of various known results found in \cite[Section 3.983] {GR}. Of particular interest is the case $a=0$, which reproduces \eqref{X1}. Setting $q=j-1/2, j=1,2,\dots$ we obtain 

\begin{equation}
\mapleinline{inert}{2d}{Int((cos(a*t)+1)/cosh(t),t = 0 .. infinity) =
Pi*(cosh(Pi*a)+1)/(3*cosh(1/2*Pi*a)+cosh(3/2*Pi*a))+1/2*Pi;}{%
\[
{\displaystyle \int _{0}^{\infty }} {\displaystyle \frac {
\mathrm{cos}(a\,t) + 1}{\mathrm{cosh}(t)}} \,dt={\displaystyle 
\frac {\pi \,(\mathrm{cosh}(\pi \,a) + 1)}{3\,\mathrm{cosh}(
{\displaystyle \frac {\pi \,a}{2}} ) + \mathrm{cosh}(
{\displaystyle \frac {3\,\pi \,a}{2}} )}}  + {\displaystyle 
\frac {\pi }{2}} 
\]
}
\end{equation}

equivalent to known results \cite[Eqs. 3.511(1) and 3.981(3)] {GR}, after redefining $a\rightarrow \pi/2/(2j-1)a$. In the case $q=j$, we reinterpret \eqref{var1} as a principal value integral, and by evaluating the residue at $x=0$, eventually arrive at 
\begin{equation}
\mapleinline{inert}{2d}{Int(sin(a*t)/sinh(t), t = 0 .. infinity) = (1/2)*Pi*sinh(Pi*a)/(cosh(Pi*a)+1)}{\[\displaystyle \int _{0}^{\infty }\!{\frac {\sin \left( at \right) }{\sinh \left( t \right) }}{dt}=(\frac {\pi}{2})\,{\frac { \,\sinh \left( \pi \,a \right) }{\cosh \left( \pi \,a \right) +1}}\]}
\end{equation}

after again redefining the variable $a$. This is also a known result \cite [Eq. 3.981(1) ] {GR}.
\newline
\subsection{Example: $F(z)=1/(z-ib)^s$}
The case \eqref{q=0} appears to have potential to lead to interesting results involving infinite series, so, setting 
\begin{equation}
F(z)=1/(z-ib)^s
\label{cut-pole}
\end{equation}
 with $b>0$ and equating real and imaginary parts eventually leads to

\begin{equation}
\mapleinline{inert}{2d}{2*cosh(Pi*r)/Pi*Int(cos(s*arctan((a*x-b)/a/x^2))*sin(2*Pi*r*x)/((a^2*
x^4+a^2*x^2-2*a*x*b+b^2)^(1/2*s))/(cosh(2*Pi*r)-cos(4*Pi*r*x)),x =
-infinity ..
infinity)+2*sinh(Pi*r)/Pi*Int(-cos(2*Pi*r*x)*sin(s*arctan((a*x-b)/a/x^
2))/((a^2*x^4+a^2*x^2-2*a*x*b+b^2)^(1/2*s))/(cosh(2*Pi*r)-cos(4*Pi*r*x
)),x = -infinity .. infinity) =
2^(2*s)*Sum(sin(s*arctan(4*b*r^2/a/(r^2+n^2)))*(-1)^n*((a^2*r^4+2*a^2*
r^2*n^2+a^2*n^4+16*b^2*r^4)/r^4)^(-1/2*s),n = 1 ..
infinity)/Pi/r+sin(s*arctan(4*b/a))*2^(2*s-1)/((a^2+16*b^2)^(1/2*s))/P
i/r;}{%
\maplemultiline{
{\displaystyle {2\,\mathrm{cosh}(\pi \,r)}{}} \,
{\displaystyle \int _{ - \infty }^{\infty }} {\displaystyle 
\frac {\mathrm{cos}(s\,\mathrm{arctan}({\displaystyle \frac {a\,x
 - b}{a\,x^{2}}} ))\,\mathrm{sin}(2\,\pi \,r\,x)}{(a^{2}\,x^{4}
 + (a\,x-b)^{2})^{s/2}\,(\mathrm{
cosh}(2\,\pi \,r) - \mathrm{cos}(4\,\pi \,r\,x))}} \,dx \\
\mbox{} - {\displaystyle {2\,\mathrm{sinh}(\pi \,r)}{}
} \,{\displaystyle \int _{ - \infty }^{\infty }}  
{\displaystyle \frac {\mathrm{cos}(2\,\pi \,r\,x)\,\mathrm{sin}(s
\,\mathrm{arctan}({\displaystyle \frac {a\,x - b}{a\,x^{2}}} ))}{
(a^{2}\,x^{4} +(a\,x-b)^{2})^{s/2}\,(\mathrm{cosh}(2\,\pi \,r) - \mathrm{cos}(4\,\pi \,r\,x))}} 
\,dx  \\
={\displaystyle {2^{2\,s}\,r^{2\,s-1}\,}{} 
{\displaystyle \sum _{n=1}^{\infty }} \,(-1)^{n}\,\frac{\mathrm{sin}(s\,\mathrm{
arctan}({\displaystyle \frac {4\,b\,r^{2}}{a\,(r^{2} + n^{2})}} )
)}{({\displaystyle {a^{2}(r^{2} + n^{2})^{2} + 16\,b^{2}\,r^{4}}{}} )^{ { {s/2}}}}\, }  
\mbox{} + {\displaystyle \frac {\mathrm{sin}(s\,\mathrm{arctan}(
{\displaystyle \frac {4\,b}{a}} ))\,2^{(2\,s - 1)}}{(a^{2} + 16\,
b^{2})^{s/2}\,r}}  }
}
\label{X5Ra}
\end{equation}
and
\begin{equation}
\mapleinline{inert}{2d}{2*cosh(Pi*r)/Pi*Int(sin(s*arctan((a*x-b)/a/x^2))*sin(2*Pi*r*x)/((a^2*
x^4+a^2*x^2-2*a*x*b+b^2)^(1/2*s))/(cosh(2*Pi*r)-cos(4*Pi*r*x)),x =
-infinity ..
infinity)+2*sinh(Pi*r)/Pi*Int(cos(2*Pi*r*x)*cos(s*arctan((a*x-b)/a/x^2
))/((a^2*x^4+a^2*x^2-2*a*x*b+b^2)^(1/2*s))/(cosh(2*Pi*r)-cos(4*Pi*r*x)
),x = -infinity .. infinity) =
2^(2*s)*r^(2*s)*Sum(cos(s*arctan(4*b*r^2/a/(r^2+n^2)))*(-1)^n*(a^2*r^4
+2*a^2*r^2*n^2+a^2*n^4+16*b^2*r^4)^(-1/2*s),n = 1 ..
infinity)/Pi/r-1/2*cos(s*arctan(4*b/a))*2^(2*s)/((a^2+16*b^2)^(1/2*s))
/Pi/r;}{%
\maplemultiline{
{\displaystyle {2\,\mathrm{cosh}(\pi \,r)}{}} \,
{\displaystyle \int _{ - \infty }^{\infty }} {\displaystyle 
\frac {\mathrm{sin}(s\,\mathrm{arctan}({\displaystyle \frac {a\,x
 - b}{a\,x^{2}}} ))\,\mathrm{sin}(2\,\pi \,r\,x)}{(a^{2}\,x^{4}
 + (a\,x-b)^{2})^{s/2}\,(\mathrm{
cosh}(2\,\pi \,r) - \mathrm{cos}(4\,\pi \,r\,x))}} \,dx \\
\mbox{} + {\displaystyle {2\,\mathrm{sinh}(\pi \,r)}{}
} \,{\displaystyle \int _{ - \infty }^{\infty }} {\displaystyle 
\frac {\mathrm{cos}(2\,\pi \,r\,x)\,\mathrm{cos}(s\,\mathrm{
arctan}({\displaystyle \frac {a\,x - b}{a\,x^{2}}} ))}{(a^{2}\,x
^{4} +(a\,x-b)^{2})^{s/2}\,(
\mathrm{cosh}(2\,\pi \,r) - \mathrm{cos}(4\,\pi \,r\,x))}} \,dx
 \\
={\displaystyle {2^{(2\,s)}\,r^{(2\,s-1)}}{}} 
{\displaystyle \sum _{n=1}^{\infty }}(-1)^{n}\, \,\frac {\displaystyle \mathrm{cos}(
s\,\mathrm{arctan}({\displaystyle \frac {4\,b\,r^{2}}{a\,(r^{2}
 + n^{2})}} ))}{\displaystyle (a^{2}(r^{2} + n^{2})^{2} + 16\,b^{2}\,r^{4})^{s/2}}\,  \mbox{} + {\displaystyle } \,{\displaystyle \frac {
\mathrm{cos}(s\,\mathrm{arctan}({\displaystyle \frac {4\,b}{a}} )
)\,2^{(2\,s-1)}}{(a^{2} + 16\,b^{2})^{s/2}\,r}} 
 \,. }
}
\label{X5Ia}
\end{equation}

In \eqref{X5Ra}, we find that the dominant left-hand terms diverge like $1/r$ as $r\rightarrow 0$, the first term on the right diverges less quickly in the same limit if $\mathbbm{R}(s)>0$, and the second term on the right simply diverges as $1/r$. Thus, by equating leading terms that diverge to the same order, we find simpler results

\begin{equation}
\mapleinline{inert}{2d}{Int((2*cos(s*arctan((a*x-b)/a/x^2))*x-sin(s*arctan((a*x-b)/a/x^2)))/(
(a^2*x^4+a^2*x^2-2*a*x*b+b^2)^(1/2*s))/(4*x^2+1),x = -infinity ..
infinity) =
Pi*sin(s*arctan(4*b/a))*2^(2*s-1)/((a^2+16*b^2)^(1/2*s));}{%
\[
{\displaystyle \int _{ - \infty }^{\infty }} {\displaystyle 
\frac {2\,\mathrm{cos}(s\,\mathrm{arctan}({\displaystyle \frac {a
\,x - b}{a\,x^{2}}} ))\,x - \mathrm{sin}(s\,\mathrm{arctan}(
{\displaystyle \frac {a\,x - b}{a\,x^{2}}} ))}{(a^{2}\,x^{4} +( a
\,x-b)^{2})^{s/2}\,(4\,x^{2} + 1
)}} \,dx={\displaystyle \frac {\pi \,\mathrm{sin}(s\,\mathrm{
arctan}({\displaystyle \frac {4\,b}{a}} ))\,2^{(2\,s - 1)}}{(a^{2
} + 16\,b^{2})^{s/2}}} 
\]
}
\label{X5Rar0}
\end{equation}

and

\begin{equation}
\mapleinline{inert}{2d}{Int((2*sin(s*arctan((a*x-b)/a/x^2))*x+cos(s*arctan((a*x-b)/a/x^2)))/(
(a^2*x^4+a^2*x^2-2*a*x*b+b^2)^(1/2*s))/(4*x^2+1),x = -infinity ..
infinity) =
Pi*cos(s*arctan(4*b/a))*2^(2*s-1)/((a^2+16*b^2)^(1/2*s));}{%
\[
{\displaystyle \int _{ - \infty }^{\infty }} {\displaystyle 
\frac {2\,\mathrm{sin}(s\,\mathrm{arctan}({\displaystyle \frac {a
\,x - b}{a\,x^{2}}} ))\,x + \mathrm{cos}(s\,\mathrm{arctan}(
{\displaystyle \frac {a\,x - b}{a\,x^{2}}} ))}{(a^{2}\,x^{4} +( a
\,x-b)^{2})^{s/2}\,(4\,x^{2} + 1
)}} \,dx={\displaystyle \frac {\pi \,\mathrm{cos}(s\,\mathrm{
arctan}({\displaystyle \frac {4\,b}{a}} ))\,2^{(2\,s - 1)}}{(a^{2
} + 16\,b^{2})^{s/2}} \,\,,\mathfrak{Re}(s)>-\frac{1}{2}\,.} 
\]
}
\label{X5Iar0}
\end{equation} 

Further results can be obtained by taking various limits of other parameters. In \eqref{X5Iar0} equate terms on both sides that vanish as $s$ when $\mathfrak{Re}(s)\rightarrow 0$, to obtain

\begin{equation}
\mapleinline{inert}{2d}{2/Pi*Int((arctan((a*x-b)/a/x^2)*x-1/4*ln(a^2*x^4+a^2*x^2-2*a*x*b+b^2)
)/(4*x^2+1),x = -infinity .. infinity) = ln(2)-1/4*ln(a^2+16*b^2);}{%
\[
{\displaystyle \frac {2}{\pi }} \,{\displaystyle \int _{ - \infty
 }^{\infty }} {\displaystyle \frac {x\,\mathrm{arctan}(
{\displaystyle \frac {a\,x - b}{a\,x^{2}}} ) - {\displaystyle 
\frac {1}{4}} \,\mathrm{ln}(a^{2}\,x^{4} +(a\,x-b)^{2})}{4\,x^{2} + 1}} \,dx=\mathrm{ln}(2) - 
{\displaystyle \frac {1}{4}} \,\mathrm{ln}(a^{2} + 16\,b^{2})
\]
}
\label{S0I}
\end{equation}

valid for $\{a,b\} \in \mathbbm{C}, a^2+16b^2 >0 $. Further results follow by similar manipulations. Comparing the second order terms $(s^2)$ in \eqref{X5Iar0} leads to a fairly complicated result, but in the limit $ b\rightarrow 0$ we find

\begin{equation}
\mapleinline{inert}{2d}{Int(ln(a^2*(v^2+1)*v^2)*(8*v*arctan(1/v)-ln(a^2*(v^2+1)*v^2))/(4*v^2+
1),v = 0 .. infinity) = -1/12*Pi^3-Pi*(ln(4/a)^2+LI[2](1/3));}{%
\[
{\displaystyle \int _{0}^{\infty }} {\displaystyle \frac {
\mathrm{ln}(a^{2}\,(v^{2} + 1)\,v^{2})\,(8\,v\,\mathrm{arctan}(
{\displaystyle1/v} ) - \mathrm{ln}(a^{2}\,(v^{2} + 1)\,
v^{2}))}{4\,v^{2} + 1}} \,dv= - {\displaystyle \frac {\pi ^{3}}{
12}}  - \pi \,(\mathrm{ln}^{2}({\displaystyle 4/a} ) + {
\mathit{Li}_{2}}({\displaystyle 1/3} ))
\]
}
\label{S2sa}
\end{equation}

valid for $\mathfrak{Re}(a)>0$. Or, instead of expanding about $s=0$, expanding about $s=1$, corresponding to a simple pole in \eqref{cut-pole}, leads to

\begin{equation}
\mapleinline{inert}{2d}{1/Pi^2*Int((3*a*x-2*b)*x/(a^2*x^4+a^2*x^2-2*a*x*b+b^2)/(4*x^2+1),x =
-infinity .. infinity) = 2*a/(a^2+16*b^2)/Pi;}{%
\[
{\displaystyle \int _{ - 
\infty }^{\infty }} {\displaystyle \frac {(3\,a\,x - 2\,b)\,x}{(a
^{2}\,x^{4} + (a\,x-b)^{2})\,(4\,x^{2} + 1)}
} \,dx={\displaystyle \frac {2\,a\,\pi}{(a^{2} + 16\,b^{2})}\,\,} 
\]
}
\label{S2a}
\end{equation}
which, Maple \cite{Maple} alternatively finds, can be expressed as the sum of the roots of a fourth order transcendental equation. 
\newline

Alternatively, by evaluating the limit $b \rightarrow 0$ in \eqref{X5Ia} we find

\begin{equation}
\mapleinline{inert}{2d}{4*cosh(Pi*r)*Int(v^(-s)*(v^2+1)^(-1/2*s)*sin(s*arctan(1/v))*sin(2*Pi*
r*v)/(-cosh(2*Pi*r)+cos(4*Pi*r*v)),v = 0 ..
infinity)+4*sinh(Pi*r)*Int(cos(2*Pi*r*v)*cos(s*arctan(1/v))*v^(-s)*(v^
2+1)^(-1/2*s)/(-cosh(2*Pi*r)+cos(4*Pi*r*v)),v = 0 .. infinity) =
-2^(2*s)*r^(2*s-1)*Sum((-1)^n*(r^2+n^2)^(-s),n = 1 ..
infinity)-2^(2*s-1)/r;}{%
\maplemultiline{
4\,\mathrm{cosh}(\pi \,r)\,{\displaystyle \int _{0}^{\infty }} 
{\displaystyle \frac {v^{ - s}\,(v^{2} + 1)^{( - s/2)}
\,\mathrm{sin}(s\,\mathrm{arctan}({\displaystyle 1/v} ))
\,\mathrm{sin}(2\,\pi \,r\,v)}{\mathrm{cos}(4\,\pi \,r\,v) - \mathrm{cosh}(2\,\pi \,r) 
}} \,dv \\
\mbox{} + 4\,\mathrm{sinh}(\pi \,r)\,{\displaystyle \int _{0}^{
\infty }} {\displaystyle \frac {\mathrm{cos}(2\,\pi \,r\,v)\,
\mathrm{cos}(s\,\mathrm{arctan}({\displaystyle 1/v} ))\,
v^{ - s}\,(v^{2} + 1)^{( - s/2)}}{ \mathrm{cos}(4\,\pi \,r\,v) - \mathrm{cosh}(2\,
\pi \,r) }} \,dv \\
= 2^{(2\,s)}\,r^{(2\,s - 1)}\,{\displaystyle \sum _{n=1}^{
\infty }} \,(-1)^{n+1}\,(r^{2} + n^{2})^{( {- s})} - {\displaystyle 
\frac {2^{(2\,s - 1)}}{r}}  }
\label{X4IB2}
}
\end{equation}

where convergence requires $1/2 <s < 1$ on the left and right hand sides, respectively, so each provides the analytic continuation of the other outside the shared region of equality. In the limit $ r \rightarrow 0$, by equating the coefficients of the terms in $1/r$ with the criterion that $s>1/2$ we find

\begin{equation}
\mapleinline{inert}{2d}{Int((v^2+1)^(-1/2*s)*(cos(s*arctan(1/v))*v^(-s)+2*v^(-s+1)*sin(s*arct
an(1/v)))/(1+4*v^2),v = 0 .. infinity) = 2^(2*s-2)*Pi;}{%
\[
{\displaystyle \int _{0}^{\infty }} {\displaystyle \frac {\mathrm{cos}(s\,\mathrm{arctan}(
{\displaystyle 1/v} )) + 2\,v \,
\mathrm{sin}(s\,\mathrm{arctan}({\displaystyle1/v} ))
}{v^{s}\, (v^{2} + 1)^{ s/2} (1 + 4\,v^{2})}} \,dv=2^{(2\,s - 2)}\,\pi 
\]
\label{Arm1}
}
\end{equation}

which is numerically valid for $-1/2<s<1$, and extendible to other values of $s$ by the principle of analytic continuation. This result may profitably be compared to \eqref{X2a} and \eqref{eta}, where the same (limiting) function $F(z)$ has been used in a different variation of the Master equation.
\newline
\subsection{Example: $F(z)=z^{s}\sin{z}$ and $F(z)=z^{s}\cos{z}$ }
Using a more complicated form of the function $F(z)$ leads to other interesting results, that are believed to be new. In \eqref{A1} set $F(z)=z^{s}\sin{z}$ with $\mathfrak{Re}(s)>-1, a \leq 2k\pi$ yielding

\begin{equation}
\mapleinline{inert}{2d}{SI2 :=
Int(v^s*(v^2+1)^(1/2*s)*(sin(s*arctan(1/v))*sin(a*v^2)*cosh(a*v)+cos(s
*arctan(1/v))*cos(a*v^2)*sinh(a*v))/sinh(2*k*Pi*v),v = 0 .. infinity)
=
1/4*(-1)^(k+1)*4^(-s)*sin(1/4*a)/k+1/2*4^(-s)*(-1)^(k+1)*k^(-2*s-1)*Sum((k^2-
n^2)^s*sin(1/4*a*(k^2-n^2)/k^2)*(-1)^n,n = 1 .. k-1);}{%
\maplemultiline{
{\displaystyle \int _{0}^{\infty }} 
{\displaystyle \frac {v^{s}\,(v^{2} + 1)^{s/2}\,(
\mathrm{sin}(s\,\mathrm{arctan}({\displaystyle 1/v } ))\,
\mathrm{sin}(a\,v^{2})\,\mathrm{cosh}(a\,v) + \mathrm{cos}(s\,
\mathrm{arctan}({\displaystyle1/v } ))\,\mathrm{cos}(a\,
v^{2})\,\mathrm{sinh}(a\,v))}{\mathrm{sinh}(2\,k\,\pi \,v)}} 
dv \\
=  {\displaystyle \frac {(-1)^{(k+1)}
\,\mathrm{sin}({\displaystyle a/4} )}{4^{s+1}k}} 
 + {\displaystyle \frac {(-1)^{(k+1)}}{(2\,k)^{1+2\,s}}}\, {\displaystyle \sum _{n=1}^{k - 1}} \,(-1)^{n}\,(k^{
2} - n^{2})^{s}\,\mathrm{sin}({\displaystyle \frac {a\,(k^{2} - n
^{2})}{4\,k^{2}}} )   }
}
\label{SI2}
\end{equation}

and, with $s=0$ and a rescaling of variables we find

\begin{equation}
\mapleinline{inert}{2d}{R2 := Int(cos(a*v^2)*sinh(2*a*v)/sinh(k*Pi*v),v = 0 .. infinity) =
(-1)^(k+1)*Sum((-1)^n*sin(a*(k+n)*(k-n)/k^2),n = 1 ..
k-1)/k-1/2*(-1)^k*sin(a)/k;}{%
\maplemultiline{
{\displaystyle \int _{0}^{\infty }} 
{\displaystyle \frac {\mathrm{cos}(a\,v^{2})\,\mathrm{sinh}(2\,a
\,v)}{\mathrm{sinh}(k\,\pi \,v)}} \,dv= 
 {\displaystyle \frac {(-1)^{k+1}}{k} \, {\displaystyle 
\sum _{n=1}^{k - 1}} \,(-1)^{n}\,\mathrm{sin}({\displaystyle 
\frac {a\,(k^2 - n^2)}{k^{2}}} )  }  + 
{\displaystyle } \,{\displaystyle \frac {(-1)^{k+1}\,
\mathrm{sin}(a)}{2\,k}}  }
\label{A1R1}
}
\end{equation}

if  $a\in\mathbbm{R}$, together with the rescaled caveat that, when $a>k\pi/2$ the left-hand side diverges. Further interesting results can be obtained by evaluating the limit $a=k\pi/2$ - see \cite{G&M}. To first order in the variable $a$ in the limit $a \rightarrow 0$ in \eqref{SI2} we find for $\Re(s)>-1$

\begin{equation}
\mapleinline{inert}{2d}{SI2a :=
sin(1/2*s*Pi)*Int((v^2+1)^(1/2*s)*v^(s+1)*(cos(s*arctan(v))*v+sin(s*ar
ctan(v)))/sinh(2*k*Pi*v),v = 0 ..
infinity)+cos(1/2*s*Pi)*Int((v^2+1)^(1/2*s)*v^(s+1)*(-sin(s*arctan(v))
*v+cos(s*arctan(v)))/sinh(2*k*Pi*v),v = 0 .. infinity) =
-1/8*4^(-s)*(-1)^k*k^(-2*s-3)*Sum((-1)^n*(k^2-n^2)^(s+1),n = 1 ..
k-1)-(-1)^k*4^(-s-2)/k;}{%
\maplemultiline{
\mathrm{sin}({\displaystyle \frac {s\,\pi }{2}} 
)\,{\displaystyle \int _{0}^{\infty }} {\displaystyle \frac {(v^{
2} + 1)^{s/2}\,v^{(s + 1)}\,(v\,\mathrm{cos}(s\,\mathrm{
arctan}(v)) + \mathrm{sin}(s\,\mathrm{arctan}(v)))}{\mathrm{
sinh}(2\,k\,\pi \,v)}} \,dv \\
\mbox{} + \mathrm{cos}({\displaystyle \frac {s\,\pi }{2}} )\,
{\displaystyle \int _{0}^{\infty }} {\displaystyle \frac {(v^{2}
 + 1)^{s/2}\,v^{(s + 1)}\,( -v\, \mathrm{sin}(s\,\mathrm{
arctan}(v)) + \mathrm{cos}(s\,\mathrm{arctan}(v)))}{\mathrm{
sinh}(2\,k\,\pi \,v)}} \,dv \\
= {\displaystyle } {\displaystyle \frac{(-1)^{(k+1)}}{(2\,k)^{2\,s+3}}}\,{\displaystyle \sum _{n=1}^{k - 1}} \,(-1)^{n}\,(k^{2
} - n^{2})^{s + 1} - {\displaystyle \frac {(-1)^{k}}{2^{2\,s+4}k}}  }
}
\label{SI2b}
\end{equation}
\newline

Alternatively, in \eqref{A1} set $F(z)=z^{s}\cos{z}$ with $\mathfrak{Re}(s)>0$ and $a<2k\pi$ yielding (see also \eqref{SI2})
\begin{equation}
\mapleinline{inert}{2d}{CI2 :=
Int(v^s*(v^2+1)^(1/2*s)*(sin(s*arctan(1/v))*cos(a*v^2)*cosh(a*v)-cos(s
*arctan(1/v))*sin(a*v^2)*sinh(a*v))/sinh(2*k*Pi*v),v = 0 .. infinity)
=
-1/4*(-1)^k*4^(-s)*cos(1/4*a)/k-1/2*(-1)^k*k^(-2*s-1)*4^(-s)*Sum((k^2-
n^2)^s*(-1)^n*cos(1/4*a*(k^2-n^2)/k^2),n = 1 .. k-1);}{%
\maplemultiline{
{\displaystyle \int _{0}^{\infty }} 
{\displaystyle \frac {v^{s}\,(v^{2} + 1)^{s/2}\,(
\mathrm{sin}(s\,\mathrm{arctan}({\displaystyle 1/v} ))\,
\mathrm{cos}(a\,v^{2})\,\mathrm{cosh}(a\,v) - \mathrm{cos}(s\,
\mathrm{arctan}({\displaystyle 1/v} ))\,\mathrm{sin}(a\,
v^{2})\,\mathrm{sinh}(a\,v))}{\mathrm{sinh}(2\,k\,\pi \,v)}}  
dv\\= {\displaystyle \frac {(-1)^{
(k+1)}\,\mathrm{cos}({\displaystyle a/4} )}{2^{2\,s+2}k}} 
\mbox{} + {\displaystyle \frac {(-1)^{(k+1)}}{(2\,k)^{2\,s+1}}} \,{\displaystyle \sum _{n=1}^{k - 1}} 
\,(-1)^{n}\,(k^{2} - n^{2})^{s}\,\mathrm{cos}({\displaystyle 
\frac {a\,(k^{2} - n^{2})}{4\,k^{2}}} ) .}
}
\label{CI2}
\end{equation}
Notice that the right-hand side of \eqref{CI2}\ provides the analytic continuation of the left-hand side when $\mathfrak{Re}(s) \leq 0$. In particular, for the case $s=0$, simple substitution $s=0$ into the left-hand side of \eqref{CI2}, although tempting, cannot yield a valid integral representation for this case since the integral diverges. In the case $a=0,\,   \mathfrak{Re}(s)>0$, we find
\begin{equation}
\mapleinline{inert}{2d}{SI2b :=
sin(1/2*s*Pi)*Int((t^2+1)^(1/2*s)*cos(s*arctan(t))/(t^(-s))/sinh(2*k*P
i*t),t = 0 ..
infinity)-cos(1/2*s*Pi)*Int((t^2+1)^(1/2*s)*sin(s*arctan(t))/(t^(-s))/
sinh(2*k*Pi*t),t = 0 .. infinity) =
-(-1)^k*(2*k)^(-2*s-1)*Sum((k^2-n^2)^s*(-1)^n,n = 1 ..
k-1)-(-1)^k*2^(-2*s-2)/k;}{%
\maplemultiline{
\mathrm{sin}({\displaystyle \frac {s\,\pi }{2}} 
)\,{\displaystyle \int _{0}^{\infty }} {\displaystyle \frac {v^{s}\,(v^{
2} + 1)^{s/2}\,\mathrm{cos}(s\,\mathrm{arctan}(v))}{\mathrm{sinh}(2\,k\,\pi \,v)}} \,dv 
\mbox{} - \mathrm{cos}({\displaystyle \frac {s\,\pi }{2}} )\,
{\displaystyle \int _{0}^{\infty }} {\displaystyle \frac {v^{s}\,(v^{2}
 + 1)^{s/2}\,\mathrm{sin}(s\,\mathrm{arctan}(v))}{\mathrm{sinh}(2\,k\,\pi \,v)}} \,dv \\
= (-1)^{k+1}\,(2\,k)^{( - 2\,s - 1)}\,{\displaystyle \sum _{n=1}
^{k - 1}} \,(k^{2} - n^{2})^{s}\,(-1)^{n} - {\displaystyle 
\frac {(-1)^{k}\,2^{( - 2\,s - 2)}}{k}}\,\, .}
}
\label{CI2b}
\end{equation} 
\newline

As a companion to the above, consider the case $F(z)=z^{s}\sin(z)$ with \eqref{A.3} which yields

\begin{equation}
\mapleinline{inert}{2d}{X6A :=
2^(2*s)*Int(v^s*(v^2+1)^(1/2*s)*(sin(4*a*v^2)*cos(s*arctan(1/v))*cosh(
4*a*v)-cos(4*a*v^2)*sin(s*arctan(1/v))*sinh(4*a*v))/cosh(Pi*(2*k-1)*v)
,v = 0 .. infinity) =
-(-1)^k*Sum(((2*k-2*n-1)*(2*k+2*n-1)/(2*k-1)^2)^s*(-1)^n*sin(a*(2*k-2*
n-1)*(2*k+2*n-1)/(2*k-1)^2),n = 1 ..
k-1)/(2*k-1)-1/2*(-1)^k*sin(a)/(2*k-1);}{%
\maplemultiline{
2^{2\,s}{\displaystyle \int _{0}^{\infty }} 
{\displaystyle \frac {v^{s}\,(v^{2} + 1)^{s/2}\,(
\mathrm{sin}(a\,v^{2})\,\mathrm{cos}(s\,\mathrm{arctan}(
{\displaystyle 1/v} ))\,\mathrm{cosh}(a\,v) - 
\mathrm{cos}(a\,v^{2})\,\mathrm{sin}(s\,\mathrm{arctan}(
{\displaystyle 1/v} ))\,\mathrm{sinh}(a\,v))}{
\mathrm{cosh}(\pi \,(2\,k - 1)\,v)}} dv 
\\ = {\displaystyle {\frac{(-1)^{k+1}}{(2\,k - 1)^{(2\,s+1)}} \, {\displaystyle 
\sum _{n=1}^{k - 1}} \,(-1)^{n}\,{\displaystyle {((2\,k-1)^2 - (2\,n)^2)\,
}{}}^{ {s}}\,\mathrm{sin}
({\displaystyle \frac {a\,((2\,k-1)^2 - (2\,n)^2)}{4\,(
2\,k - 1)^{2}}} ) }{}}  
\mbox{}+ {\displaystyle \frac {1}{2}} \,{\displaystyle \frac {(
-1)^{k+1}\,\mathrm{sin}(a/4)}{2\,k - 1}}  }
}
\label{X6A}
\end{equation}
where $4a \leq (2 k-1) \pi$. For the case $a=(2 k-1)/(4 \pi)$ see \cite{G&M}.  
\newline

Similarly, in \eqref{Fcosh}, using $F(z)=z^{s}\cos(z)$ we find the companion result to \eqref{X6A}

\begin{equation}
\mapleinline{inert}{2d}{Y6A :=
2^(2*s)*Int((v^2*(v^2+1))^(1/2*s)*(cos(a*v^2)*cos(s*arctan(1/v))*cosh(
a*v)+sin(a*v^2)*sin(s*arctan(1/v))*sinh(a*v))/cosh(Pi*(2*k-1)*v),v = 0
.. infinity) =
-(-1)^k*Sum((-1)^n*((2*k-2*n-1)*(2*k+2*n-1)/(2*k-1)^2)^s*cos(1/4*a*(2*
k-2*n-1)*(2*k+2*n-1)/(2*k-1)^2),n = 1 ..
k-1)/(2*k-1)-1/2*(-1)^k*cos(1/4*a)/(2*k-1);}{%
\maplemultiline{
2^{(2\,s)}{\displaystyle \int _{0}^{\infty }} 
{\displaystyle \frac {v^{s}(v^{2} + 1)^{s/2}\,(
\mathrm{cos}(a\,v^{2})\,\mathrm{cos}(s\,\mathrm{arctan}(
{\displaystyle 1/v} ))\,\mathrm{cosh}(a\,v) + \mathrm{
sin}(a\,v^{2})\,\mathrm{sin}(s\,\mathrm{arctan}({\displaystyle 
1/v} ))\,\mathrm{sinh}(a\,v))}{\mathrm{cosh}(\pi \,(2\,k
 - 1)\,v)}}  
dv \\ = {\displaystyle \frac {(-1)^{k+1}}{(2\,k - 1)^{(2\,s+1)}}\, {\displaystyle 
\sum _{n=1}^{k - 1}} \,(-1)^{n}\,{\displaystyle {((2\,k -1)^2-( 2
\,n)^2)}{}} ^{s}\,\mathrm{cos}
({\displaystyle \frac {a\,((2\,k-1)^2 -( 2\,n)^2)}{4
\,(2\,k - 1)^{2}}} ) }{}  
\mbox{} + {\displaystyle \frac {(
-1)^{k+1} }{2}} \,{\displaystyle \frac {\mathrm{cos}({\displaystyle a/4} )}{2\,k - 1}} \,\,. 
 }
}
\label{Y6A}
\end{equation}

In the case $s=0$, uniting \eqref{X6A} and \eqref{Y6A}, yields
\begin{equation}
\mapleinline{inert}{2d}{Int(exp(i*a*x^2)*cosh(a*x)/cosh((2*k-1)*Pi*x),x = 0 .. infinity) =
(-1)^(k+1)*exp(1/4*i*a)*(1/2+Sum((-1)^n*exp(-i*a*n^2/(2*k-1)^2),n = 1
.. k-1))/(2*k-1);}{%
\[
{\displaystyle \int _{0}^{\infty }} {\displaystyle \frac {e^{i\,
a\,v^{2}}\,\mathrm{cosh}(a\,v)}{\mathrm{cosh}((2\,k - 1)\,\pi \,
v)}} \,dv={\displaystyle \frac{(-1)^{(k + 1)}\,\exp({{i\,a/4}{
})}}{2\,k - 1}{\displaystyle \left(  \! {\displaystyle \frac {1}{2}}  +  
{\displaystyle \sum _{n=1}^{k - 1}} \,(-1)^{n}\,\exp{( - \frac {i\,
a\,n^{2}}{(2\,k - 1)^{2}})}   \!  \right) }{}
} 
\]
}
\label{XY6B}
\end{equation}

where $\mathfrak{Im}(a)\geq 0$. An interesting special case arises from \eqref{XY6B} after setting $k=1$ and replacing $a\rightarrow i a$ in the real part, yielding
\begin{equation}
\mapleinline{inert}{2d}{Int(exp(-a*x^2)*cos(a*x)/cosh(Pi*x),x = 0 .. infinity) =
1/2*exp(-1/4*a);}{%
\[
{\displaystyle \int _{0}^{\infty }} {\displaystyle \frac {e^{( - 
a\,v^{2})}\,\mathrm{cos}(a\,v)}{\mathrm{cosh}(\pi \,v)}} \,dv=
{\displaystyle \frac {1}{2}} \,\exp{( - a/4)} \hspace{.5in} \mathfrak{Re}(a) \geq 0
\]
}
\label{XYk1}
\end{equation}

which appears to be new. The special case $a=(2\,k-1) \pi$ in \eqref{XY6B} is treated in \cite{G&M}. 
\newline

Applying the same function $F(z)=z^{s} \sin(z) $ to \eqref{A.4} leads to a generalization of the foregoing, specifically

\begin{equation}
\mapleinline{inert}{2d}{LS :=
Int(cosh(Pi*(2*k-1)*t)*(t^2*(t^2+1))^(1/2*s)*exp(-a*t)*sin(s*arctan(1/
t)+a*t^2)/(cosh(2*Pi*b)+cosh(2*Pi*(2*k-1)*t)),t = -infinity ..
infinity) =
1/cosh(Pi*b)/((2*k-1)^(2*s+1))/(2^(2*s+2))*Sum(B[n]*(-1)^n*(sin(s*T[n]
-A[n])*exp(-X[n])-sin(s*T[n]+A[n])*exp(X[n])),n = 1 .. 2*k-1);}{%
\maplemultiline{
{\displaystyle \int _{ - \infty }^{\infty }} 
{\displaystyle \frac {\mathrm{cosh}(\pi \,(2\,k - 1)\,v)\,(v^{2}
\,(v^{2} + 1))^{s/2}\,\exp{( - a\,v)}\,\mathrm{sin}(s\,
\mathrm{arctan}({\displaystyle 1/v} ) + a\,v^{2})}{
\mathrm{cosh}(2\,\pi \,b) + \mathrm{cosh}(2\,\pi \,(2\,k - 1)\,v)
}} \,dv \\
={\displaystyle \frac {{\displaystyle \sum _{n=1}^{2\,k - 1}} \,{B
_{n}}\,(-1)^{n}\, \left ( \, \mathrm{sin}(s\,{T_{n}} - {A_{n}})\,\exp{( - {X_{
n}})} - \mathrm{sin}(s\,{T_{n}} + {A_{n}})\,\exp{{X_{n}}} \, \right) }{
\mathrm{cosh}(\pi \,b)\,(2\,k - 1)^{(2\,s + 1)}\,2^{(2\,s + 2)}}
}  }
}
\label{C7}
\end{equation} 
where

\begin{maplelatex}
\mapleinline{inert}{2d}{T[n] = -arctan(2*b*(-n+k)/(b^2-n^2+1/4+2*k*n-k));}{%
\[
{T_{n}}= - \mathrm{arctan} \left(  \! {\displaystyle \frac {2\,b
\,( - n + k)}{b^{2} - n^{2} + {\displaystyle \frac {1}{4}}  + 2\,
k\,n - k}}  \!  \right) 
\]
}
\end{maplelatex}

\begin{maplelatex}
\mapleinline{inert}{2d}{A[n] = 1/4*a*(4*b^2-4*n^2+1+8*k*n-4*k)/(2*k-1)^2;}{%
\[
{A_{n}}={\displaystyle \frac {a\,(4\,b^{2} - 4\,n^{2} + 1 + 8\,k
\,n - 4\,k)}{4\,(2\,k - 1)^{2}}} 
\]
}
\end{maplelatex}

\begin{maplelatex}
\mapleinline{inert}{2d}{B[n] =
(((2*n-1)^2+4*b^2)*(4*b^2-16*k*n+4*n^2-8*k+1+16*k^2+4*n))^(1/2*s);}{%
\[
{B_{n}}=(((2\,n - 1)^{2} + 4\,b^{2})\,(4\,b^{2} - 16\,k\,n + 4\,n
^{2} - 8\,k + 1 + 16\,k^{2} + 4\,n))^{(s/2)}
\]
}
\end{maplelatex}

\begin{maplelatex}
\mapleinline{inert}{2d}{X[n] = 2*a*b*(-n+k)/(2*k-1)^2;}{%
\[
{X_{n}}={\displaystyle \frac {2\,a\,b\,( - n + k)}{(2\,k - 1)^{2}
}.} 
\]
}
\end{maplelatex}

Simpler variations are, for the case $s=0$,
\begin{equation}
\mapleinline{inert}{2d}{C1 :=
Int(cosh(Pi*(2*k-1)*v)*sin(a*v^2)*cosh(a*v)/(cosh(2*Pi*b)+cosh(2*Pi*(2
*k-1)*v)),v = 0 .. infinity) =
-1/4*Sum((-1)^n*sin(1/4*a*(4*b^2-4*n^2+1+8*k*n-4*k)/(2*k-1)^2)*cosh(2*
a*b*(-n+k)/(2*k-1)^2),n = 1 .. 2*k-1)/(2*k-1)/cosh(Pi*b);}{%
\maplemultiline{
{\displaystyle \int _{0}^{\infty }} 
{\displaystyle \frac {\mathrm{cosh}(\pi \,(2\,k - 1)\,v)\,
\mathrm{sin}(a\,v^{2})\,\mathrm{cosh}(a\,v)}{\mathrm{cosh}(2\,\pi
 \,b) + \mathrm{cosh}(2\,\pi \,(2\,k - 1)\,v)}} \,dv= 
 - {\displaystyle \frac {1}{4\, (2\,k - 1)\,\mathrm{cosh}
(\pi \,b)}} \,{\displaystyle  {
{\displaystyle \sum _{n=1}^{2\,k - 1}} \,(-1)^{n}\,\mathrm{sin}(A_{n} )\,\mathrm{cosh}(X_{n})}{}},  }
}
\label{C7a}
\end{equation}
for the case $s=1, k=1$
\begin{equation}
\mapleinline{inert}{2d}{C2b :=
int(cosh(Pi*v)*v*(cosh(a*v)*v*sin(a*v^2)-sinh(a*v)*cos(a*v^2))/(cosh(2
*Pi*b)+cosh(2*Pi*v)),v = 0 .. infinity) =
1/16*(1+4*b^2)*sin(1/4*a*(1+4*b^2))/cosh(Pi*b);}{%
\maplemultiline{
{\displaystyle \int _{0}^{\infty }} 
{\displaystyle \frac {v\,\mathrm{cosh}(\pi \,v)\,(v\,\mathrm{cosh}(a
\,v)\,\mathrm{sin}(a\,v^{2}) - \mathrm{sinh}(a\,v)\,\mathrm{
cos}(a\,v^{2}))}{\mathrm{cosh}(2\,\pi \,b) + \mathrm{cosh}(2\,\pi
 \,v)}} \,dv= 
{\displaystyle \frac {1}{16}} \,{\displaystyle \frac {(1 + 4\,b^{
2})\,\mathrm{sin}({\displaystyle \frac {a\,(1 + 4\,b^{2})}{4}} )
}{\mathrm{cosh}(\pi \,b)}}  }
}
\label{C7_s=1,k=1}
\end{equation}

and, to first order in the variable $a$, at $a=0, s=0$
\begin{equation}
\mapleinline{inert}{2d}{C0 :=
int(v^2*(v^2-1)*cosh(Pi*(2*k-1)*v)/(cosh(2*Pi*b)+cosh(2*Pi*(2*k-1)*v))
,v = 0 .. infinity) =
1/64*(1+4*b^2)*(4*b^2+1-16*k^2+16*k)/(2*k-1)^5/cosh(Pi*b);}{%
\[
{\displaystyle \int _{0}^{\infty }} 
{\displaystyle \frac {v^{2}\,(v^{2} - 1)\,\mathrm{cosh}(\pi \,(2
\,k - 1)\,v)}{\mathrm{cosh}(2\,\pi \,b) + \mathrm{cosh}(2\,\pi \,
(2\,k - 1)\,v)}} \,dv={\displaystyle \frac {1}{64}} \,
{\displaystyle \frac {(1 + 4\,b^{2})\,(4\,b^{2} + 1 - 16\,k^{2}
 + 16\,k)}{(2\,k - 1)^{5}\,\mathrm{cosh}(\pi \,b)}.} 
\]
}
\label{C7a=0}
\end{equation}

(It can be shown, independently of the above, that the generalization of \eqref{C7a=0} obtained by redefining $k \rightarrow \beta$ where $\beta \in \mathbbm{C}$  is true.) The result corresponding to \eqref{C7}, using the same identifications of variables for the case $F(z)=z^s\cos(z)$ is

\begin{equation}
\mapleinline{inert}{2d}{SEq3 :=
Int(exp(-a*t)*cosh(Pi*(2*k-1)*t)*cos(s*arctan(1/t)+a*t^2)*(t^2*(t^2+1)
)^(1/2*s)/(cosh(2*Pi*b)+cosh(2*Pi*(2*k-1)*t)),t = -infinity ..
infinity) =
-Sum(B[n]*(-1)^n*((cos(s*T[n]-A[n])+cos(s*T[n]+A[n]))*cosh(X[n])+(-cos
(s*T[n]-A[n])+cos(s*T[n]+A[n]))*sinh(X[n])),n = 1 ..
2*k-1)/(2^(2*s+2))/(2*k-1)/cosh(Pi*b);}{%
\maplemultiline{
{\displaystyle \int _{ - \infty }^{\infty }} 
{\displaystyle \frac {\exp{( - a\,v)}\,\mathrm{cosh}(\pi \,(2\,k - 
1)\,v)\,\mathrm{cos}(s\,\mathrm{arctan}({\displaystyle 1/v} ) + a\,v^{2})\,(v^{2}\,(v^{2} + 1))^{s/2}}{
\mathrm{cosh}(2\,\pi \,b) + \mathrm{cosh}(2\,\pi \,(2\,k - 1)\,v)
}} \,dv\\={ \displaystyle \frac { -1}{ 2^{(2\,s + 1)}\,(2\,k - 1)^{2\,s+1}\,\mathrm{cosh}(\pi \,b)}} \\ \times {\displaystyle \sum _{n=1}^{2\,k - 1}} {B_{n}}\,(-1)^{
n}  
(\mathrm{cos}(s\,{T_{n}} ) \, \mathrm{cos}(A_{n})\,\mathrm{cosh}({X_{n}})   - \mathrm{sin}(s\,{T_{n}}
) \, \mathrm{sin}({A_{n}})\,\mathrm{sinh}(
{X_{n}}) ). }
}
\label{C7A}
\end{equation}

Since \eqref{C7A} corresponds to \eqref{C7} with the replacement $ sin \rightarrow cos$ throughout, the right-hand side has been rearranged to display an alternate structure. A new variant of \eqref{C7a} is obtained when $s=0$, yielding
\begin{equation}
\mapleinline{inert}{2d}{C7A2a :=
Int(cosh(Pi*(2*k-1)*v)*cos(a*v^2)*cosh(a*v)/(cosh(2*Pi*b)+cosh(4*Pi*v*
k-2*Pi*v)),v = 0 .. infinity) =
-1/4*Sum((-1)^n*cosh(2*a*b*(-n+k)/(2*k-1)^2)*cos(1/4*a*(4*b^2-4*n^2+1+
8*k*n-4*k)/(2*k-1)^2),n = 1 .. 2*k-1)/(2*k-1)/cosh(Pi*b);}{%
\maplemultiline{
{\displaystyle \int _{0}^{\infty }} 
{\displaystyle \frac {\mathrm{cosh}(\pi \,(2\,k - 1)\,v)\,
\mathrm{cos}(a\,v^{2})\,\mathrm{cosh}(a\,v)}{\mathrm{cosh}(2\,\pi
 \,b) + \mathrm{cosh}(2\,\pi\,v\,(2\,k-1))}} \,dv= 
 - {\displaystyle \frac {1}{4\,(2\,k - 1)\,\mathrm{cosh}
(\pi \,b)}} \,{\displaystyle {
{\displaystyle \sum _{n=1}^{2\,k - 1}} \,(-1)^{n}\,\mathrm{cosh}(X_{n} )\,
\mathrm{cos}(A_{n} )}\,.{}}  }
}
\label{C7A2}
\end{equation}
For the case $b=0,a=(2k-1)\pi$ see \cite{G&M}.

\subsection{Example: $F(z)=z^s \exp(-bz)$}

Setting $F(z)=z^s \exp(-b z)$ in \eqref{gen2} leads to a large number of possibilities. The general result is fairly lengthy, but special cases lead to results that appear to be new. Without loss of generality, let $a=1$ since it always appears in the form $a b$. From the real part, with $q=0$ we find
\begin{equation}
\mapleinline{inert}{2d}{C5RB :=
Int(exp(-b*v^2)*v^s*(v^2+1)^(1/2*s)*(-cos(b*v-s*arctan(1/v))*cos(2*Pi*
r*v)*sinh(Pi*r)-sin(b*v-s*arctan(1/v))*sin(2*Pi*r*v)*cosh(Pi*r))/(-cos
(2*Pi*r*v)^2+cosh(Pi*r)^2),v = 0 .. infinity) =
-2^(-2*s-2)*exp(-1/4*b)/r-1/2*exp(-1/4*b)*Sum((-1)^n*(n^2+r^2)^s*exp(-
1/4*b*n^2/r^2),n = 1 .. infinity)/((2*r)^(2*s))/r;}{%
\maplemultiline{
{\displaystyle \int _{0}^{\infty }} \exp{( - b\,v
^{2})}\,v^{s}\,(v^{2} + 1)^{s/2}( \mathrm{cos}(b\,v
 - s\,\mathrm{arctan}({\displaystyle 1/v} ))\,\mathrm{
cos}(2\,\pi \,r\,v)\,\mathrm{sinh}(\pi \,r) \\
\mbox{} + \mathrm{sin}(b\,v - s\,\mathrm{arctan}({\displaystyle 
1/v} ))\,\mathrm{sin}(2\,\pi \,r\,v)\,\mathrm{cosh}(\pi 
\,r)) \left/ {\vrule height0.44em width0em depth0.44em}
 \right. \!  \! ( \mathrm{
cosh^{2}}(\pi \,r) - \mathrm{cos^{2}}(2\,\pi \,r\,v) )dv \\
= {\displaystyle \frac {2^{( - 2\,s - 2)}\,\exp{( - b/4)}
}{r}}  + {\displaystyle \frac {\exp{
( - b/4)}}{(2\,r)^{2s+1}}} \,{\displaystyle 
{\, {\displaystyle \sum _{n=1}^{
\infty }} \,(-1)^{n}\,(n^{2} + r^{2})^{s}\,\exp{( - \frac {b\,n^{2}
}{4\,r^{2}})}   }{}}\, ,  }
}
\label{C00}
\end{equation}

in the special case $q=0, s=0$, we find
\begin{equation}
\mapleinline{inert}{2d}{C0 :=
Int(exp(-b*v^2)*(-sinh(Pi*r)*cos(2*Pi*r*v)*cos(b*v)-cosh(Pi*r)*sin(2*P
i*r*v)*sin(b*v))/(-cos(2*Pi*r*v)^2+cosh(Pi*r)^2),v = 0 .. infinity) =
-1/4*exp(-1/4*b)/r-1/2*exp(-1/4*b)*Sum((-1)^n*exp(-1/4*b*n^2/r^2),n =
1 .. infinity)/r;}{%
\maplemultiline{
{\displaystyle \int _{0}^{\infty }} 
{\displaystyle \frac {\exp{( - b\,v^{2})}\,(  \mathrm{sinh}(\pi \,
r)\,\mathrm{cos}(2\,\pi \,r\,v)\,\mathrm{cos}(b\,v) + \mathrm{
cosh}(\pi \,r)\,\mathrm{sin}(2\,\pi \,r\,v)\,\mathrm{sin}(b\,v))
}{ \mathrm{cos^{2}}(2\,\pi \,r\,v) - \mathrm{cosh^{2}}(\pi \,r)
}} \,dv \\
 = - {\displaystyle \frac {1}{4}} \,{\displaystyle \frac {e^{( - 
b/4)}}{r}}  - {\displaystyle \frac {e^{(-b/4)}}{2\,r}} \,
{\displaystyle  { 
{\displaystyle \sum _{n=1}^{\infty }} \,(-1)^{n}\,\exp{( - {b
\,n^{2}}/{4\,r^{2}})} }{}}  }
}
\label{C10}
\end{equation}
and, in the limit $r \rightarrow 0$ with $\mathfrak{Re}(b)\ge0$
\begin{equation}
\mapleinline{inert}{2d}{C1 := int(exp(-b*v^2)*(cos(b*v)+2*sin(b*v)*v)/(4*v^2+1), v = 0 .. infinity) = (1/4)*Pi*exp(-(1/4)*b)}{\[\displaystyle \int _{0}^{\infty }\!{\frac {{{ \exp}{(-b{v}^{2})}} \left( \cos \left( bv \right) +2\,\sin \left( bv \right) v \right) }{4\,{v}^{2}+1}}{dv}=(\frac {\pi} {4}) \, {{\rm \exp}({-b/4})}\]\,.}
\label{C10a}
\end{equation}  
From the imaginary part, the corresponding result is also fairly complicated, but, in the case $r=0, s=0$ we find  
\begin{equation}
\mapleinline{inert}{2d}{C3A :=
Int(exp(-b*v^2)*(sinh(2*Pi*q*v)*cos(Pi*q)*sin(b*v)+cosh(2*Pi*q*v)*sin(
Pi*q)*cos(b*v))/(-cos(Pi*q)^2+cosh(2*Pi*q*v)^2),v = 0 .. infinity) =
1/4*(exp(-1/4*b)+2*Sum((-1)^n*exp(1/4*b*(n-q)*(n+q)/q^2),n = 1 ..
floor(q)))/q;}{%
\maplemultiline{
{\displaystyle \int _{0}^{\infty }} 
{\displaystyle \frac {\exp{( - b\,v^{2})}\,(\mathrm{sinh}(2\,\pi \,
q\,v)\,\mathrm{cos}(\pi \,q)\,\mathrm{sin}(b\,v) + \mathrm{cosh}(
2\,\pi \,q\,v)\,\mathrm{sin}(\pi \,q)\,\mathrm{cos}(b\,v))}{ \mathrm{cosh^{2}}(2\,\pi \,q\,v) - 
\mathrm{cos^{2}}(\pi \,q) }} \,
dv \\
={\displaystyle (\frac {1}{4\,q}) } \left ( \! \,{\displaystyle  {\exp{( - 
b/4)} + 2\, {\displaystyle \sum _{n=1}^{
\lfloor q \rfloor}} \,(-1)^{n}\,\exp{(\frac {b\,(n^2 - q^2)}{4
\,q^{2}})} }{} }  \right ) }  
}
\label{C10b}
\end{equation}

where the right-hand side represents the left-hand side if $ q = j$ and $b \in \mathbbm{C}$ is permitted. If  $ r=0 $ and $q=j, j \neq 0$ we find (with $ \mathfrak{Re}(s) > 0$ )
\begin{equation}
\mapleinline{inert}{2d}{C5 :=
Int(v^s*(v^2+1)^(1/2*s)*exp(-b*v^2)*(sin(s*arctan(1/v))*cos(b*v)-cos(s
*arctan(1/v))*sin(b*v))/sinh(2*Pi*j*v),v = 0 .. infinity) =
-(-1)^j*4^(-1-s)*(exp(-1/4*b)+2*j^(-2*s)*Sum((-1)^n*(j^2-n^2)^s*exp(-1
/4*b*(j-n)*(j+n)/j^2),n = 1 .. floor(j)))/j;}{%
\maplemultiline{
{\displaystyle \int _{0}^{\infty }} 
{\displaystyle \frac {v^{s}\,(v^{2} + 1)^{s/2}\,\exp{(
 - b\,v^{2})}\,(\mathrm{sin}(s\,\mathrm{arctan}({\displaystyle 
1/v} ))\,\mathrm{cos}(b\,v) - \mathrm{cos}(s\,\mathrm{
arctan}({\displaystyle1/v} ))\,\mathrm{sin}(b\,v))}{
\mathrm{sinh}(2\,\pi \,j\,v)}} \,dv \\
 ={\displaystyle  {(-1)^{j+1}\,4^{( - 1 - s)}\, \left(  \! \exp
{( - b / 4)} + 2\,j^{ - 2\,s}\, 
{\displaystyle \sum _{n=1}^{j}} \,(-1)^{n}\,(j^{2
} - n^{2})^{s}\,\exp{( - \frac {b\,(j^2 - n^2)}{4\,j^{2}})}
\!  \right) }{/j}} . }
}
\label{C10c}
\end{equation}

Note that \eqref{C10c} is discontinuous as a function of $s$ at $s=0$, and that \eqref{C10b} is the corresponding result for that case. See also \eqref{A11}. In the case $b=0$ we obtain 
\begin{equation}
\mapleinline{inert}{2d}{C5a :=
int((-1)^j*v^s*(v^2+1)^(1/2*s)*sin(s*arctan(1/v))/sinh(2*Pi*j*v),v = 0
.. infinity) = -4^(-1-s)*(1+2*j^(-2*s)*sum((-1)^n*(j^2-n^2)^s,n = 1 ..
j))/j;}{%
\maplemultiline{
{\displaystyle \int _{0}^{\infty }} 
{\displaystyle \frac {v^{s}\,(v^{2} + 1)^{s/2
}\,\mathrm{sin}(s\,\mathrm{arctan}({\displaystyle 1/v} 
))}{\mathrm{sinh}(2\,\pi \,j\,v)}} \,dv= 
 (-1)^{j+1}\, {\displaystyle {4^{( - 1 - s)}\,\left( 1 + 2\,j^{ - 2\,s}\,
{\displaystyle \sum _{n=1}^{j}} \,(-1)^{n}\,(j^{2} - n^{2})^{s}) \right)
}{/j}}  }
}
\label{C10d}
\end{equation}
which should be compared with \eqref{x2a2} and \eqref{x2b2} with reference to \eqref{arctan}.

\subsection{Example: $F(z)=\frac{\Gamma(a+z)}{\Gamma(b+z)}$}
As a more complicated illustrative example, consider
\begin{equation}
F(z)=\frac{\Gamma(a+z)}{\Gamma(b+z)}
\label{GamRat}
\end{equation}
with $b>a>1/2$ for which the only singularities are simple poles $z=-(a+k)$, $k=0,1,2,\cdots$.
From \eqref{c_Inf} one has
\begin{equation}
\int_{-1/4-i\infty}^{-1/4+i\infty}\frac{\Gamma(a+s(1-s))}{\Gamma(b+s(1-s))}\sec \pi s\frac{ds}{2\pi i}=\frac{\Gamma(a+1/4)}{2\pi\Gamma(b+1/4)}.
\label{XGam1}
\end{equation}
Since the integrand in \eqref{XGam1} decays exponentially (off the real axis) at infinity in the right half $s-$plane, the contour can be closed by a large semicircle. By summing the residues at the simple poles (all real) thus enclosed: $s=k+1/2$ due to $\cos \pi s$ and $s=[1+\sqrt{1+4(a+k)}]/2$ due to the gamma function,
$k=0,1,2,\cdots$,
one obtains
\begin{equation}
\sum_{n=0}^{\infty} (-1)^n\frac{\Gamma(a-n^2)}{\Gamma(b-n^2)}=\frac{1}{2}\frac{\Gamma(a)}{\Gamma(b)}
+\frac{\pi}{2}\sum_{n=0}^{\infty}\frac{(-1)^n}{n!}\frac{\csc[\pi\sqrt{a+n}]}{\sqrt{a+n}\Gamma(b-a-n)}
\label{XGam2}
\end{equation}
where $a$ and $b$ have been shifted by 1/4. In particular, if $a$ and $b$ differ by a positive integer $k$, one has the closed-form evaluation of the series
\begin{equation}
\sum_{n=0}^{\infty}(-1)^n\frac{\Gamma(a-n^2)}{\Gamma(a+k-n^2)}=\sum_{n=0}^{\infty}\frac{(-1)^n}{\prod_{j=0}^{k-1}(a+j-n^2)}$$
$$=\frac{1}{2}\frac{\Gamma(a)}{\Gamma(a+k)}
+\frac{\pi}{2(k-1)!}\sum_{n=0}^{k-1}(-1)^n\left(\begin{array}{c} k-1\\
n
\end{array}\right)\frac{\csc[\pi\sqrt{a+n}]}{\sqrt{a+n}},
\label{XGam3}
\end{equation}
which appears to be unlisted previously except for $k=0,1$, although both Maple \cite{Maple} and Mathematica \cite{Math} can evaluate the left-hand sum for any particular value of $k$ that was tried.  It seems clear that \eqref{XGam1}-\eqref{XGam3} can be generalized to ratios of arbitrarily many Gamma functions including further parameters. 

\section{Summary and comments}
In the second section of this paper, a new master theorem was presented for the analytic integration of an almost arbitrary function over the real line, when the integrand includes the reciprocal $cosh$ function . The principles used to obtain that result were used to generalize the theorem to include other combinations of reciprocal $cosh$ and $sinh$ kernels. In the third section, each of the theorems so-obtained was tested by employing simple ``arbitrary'' functions to obtain new evaluations of a general class of integrals, and, in several cases, a finite sum of trigonometric functions, which are dealt with elsewhere \cite{G&M}. A number of variations were presented in the second section, and, in Appendix A, another set of simple theorems and results were presented that also yield what appear to be new evaluations of interesting integrals. As shown in the last example, it is possible to use the Master theorem to transform an infinite series, and in a special case, truncate it. This appears to be possible by expanding the contour and identifying the residues. In a future work, we will turn to exploring other applications of (1.3) and its brethren to special functions, with particular emphasis on an examination of what can be learned about their zeros. 

\vskip 0.2in
\noindent
{\bf Acknowledgement}
The work of  MLG  has been partially supported by the Spanish Ministerio de Educaci\'on y Ciencia (Project MTM 2005-09183) and he wishes to thank  the University of Valladolid for continuing hospitality.

\appendix
\section{Appendix A (Variations)}

In Section 2, a number of variations were presented. The following demonstrates the utility of these variations, by, first of all, presenting the general equations for each of some obvious variations. Following, for a special choice of the function $F(z)$, a variety of results are obtained. The generalized master equations are, for an obvious set of variations:
\small
\begin{equation}
\mapleinline{inert}{2d}{Int(F(x*(i+x))/sinh(Pi*p*(2*x+i)),x = -infinity .. infinity) =
1/2*i/p*Sum((-1)^n*F(-1/4*((n+1)^2-p^2)/p^2),n = \lfloor -p \rfloor ..
-1+floor(p));}{%
\[
{\displaystyle \int _{ - \infty }^{\infty }} {\displaystyle 
\frac {\mathit{F}(x\,(i + x))}{\mathrm{sinh}(\pi \,p\,(2\,x + i))
}} \,dx={\displaystyle  {{\displaystyle \frac {i}{2\,p}} \,
 {\displaystyle \sum _{n=\mathrm{\lfloor} - p \rfloor}^{ - 1 + \mathrm{ \lfloor}p \rfloor}} \,(-1)^{n}\,\mathit{F}( - {\displaystyle 
\frac {1}{4}} \,{\displaystyle \frac {(n + 1)^{2} - p^{2}}{p^{2}}
} ) }{}} 
\]
}
\end{equation}
\begin{equation}
\mapleinline{inert}{2d}{Int(Fe(x)*Fe(I+x)/sinh(Pi*p*(2*x+I)),x = -infinity .. infinity) =
1/2*I*Sum((-1)^n*Fe(1/2*I*(n+1-p)/p)*Fe(1/2*I*(n+1+p)/p),n = floor(-p)
.. -1+floor(p))/p;}{%
\[
{\displaystyle \int _{ - \infty }^{\infty }} {\displaystyle 
\frac {\mathit{Fe}(x)\,\mathit{Fe}(i + x)}{\mathrm{sinh}(\pi \,p
\,(2\,x + i))}} \,dx={\displaystyle {{\displaystyle \frac {
i}{2\,p}} \,{\displaystyle \sum _{n=\mathrm{ \lfloor}
 - p \rfloor}^{ - 1 + \mathrm{\lfloor}p\rfloor}} \,(-1)^{n}\,\mathit{Fe}
 \left(  \! {\displaystyle \frac {{\displaystyle } \,
i\,(n + 1 - p)}{2\,p}}  \!  \right) \,\mathit{Fe} \left(  \! 
{\displaystyle \frac {{\displaystyle } \,i\,(n + 1 + 
p)}{2\,p}}  \!  \right)  }{}} 
\]
}
\end{equation}
\begin{equation}
\mapleinline{inert}{2d}{Int(Fo(x)*Fo(I+x)/sinh(Pi*p*(2*x+I)),x = -infinity .. infinity) =
1/2*I*Sum((-1)^n*Fo(1/2*I*(n+1-p)/p)*Fo(1/2*I*(n+1+p)/p),n = floor(-p)
.. -1+floor(p))/p;}{%
\[
{\displaystyle \int _{ - \infty }^{\infty }} {\displaystyle 
\frac {\mathit{Fo}(x)\,\mathit{Fo}(i + x)}{\mathrm{sinh}(\pi \,p
\,(2\,x + i))}} \,dx={\displaystyle {{\displaystyle \frac {
i}{2\,p}} \, {\displaystyle \sum _{n=\mathrm{\lfloor}
 - p \rfloor}^{ - 1 + \mathrm{\lfloor}p\rfloor}} \,(-1)^{n}\,\mathit{Fo}
 \left(  \! {\displaystyle \frac {{\displaystyle } \,
i\,(n + 1 - p)}{2\,p}}  \!  \right) \,\mathit{Fo} \left(  \! 
{\displaystyle \frac {{\displaystyle} \,i\,(n + 1 + 
p)}{2\,p}}  \!  \right)  }{}} 
\]
}
\end{equation}
\begin{equation}
\mapleinline{inert}{2d}{Int((Fe(x)+Fe(I+x))/sinh(Pi*p*(2*x+I)),x = -infinity .. infinity) =
1/2*I*Sum((-1)^n*(Fe(1/2*I*(n+1-p)/p)+Fe(1/2*I*(n+1+p)/p)),n =
floor(-p) .. -1+floor(p))/p;}{%
\[
{\displaystyle \int _{ - \infty }^{\infty }} {\displaystyle 
\frac {\mathit{Fe}(x) + \mathit{Fe}(i + x)}{\mathrm{sinh}(\pi \,p
\,(2\,x + i))}} \,dx={\displaystyle {{\displaystyle \frac {
i}{2\,p}} \, {\displaystyle \sum _{n=\mathrm{\lfloor} - p\rfloor}^{ - 1 + \mathrm{\lfloor}p\rfloor}} \,(-1)^{n}\, \left(  \! 
\mathit{Fe} \left(  \! {\displaystyle \frac {{\displaystyle 
} \,i\,(n + 1 - p)}{2\,p}}  \!  \right)  + \mathit{Fe}
 \left(  \! {\displaystyle \frac {{\displaystyle } \,
i\,(n + 1 + p)}{2\,p}}  \!  \right)  \!  \right)   }{}} 
\]
}
\end{equation}
\begin{equation}
\mapleinline{inert}{2d}{Int((Fo(x)-Fo(I+x))/sinh(Pi*p*(2*x+I)),x = -infinity .. infinity) =
1/2*I*Sum((-1)^n*(Fo(1/2*I*(n+1-p)/p)-Fo(1/2*I*(n+1+p)/p)),n =
floor(-p) .. -1+floor(p))/p;}{%
\[
{\displaystyle \int _{ - \infty }^{\infty }} {\displaystyle 
\frac {\mathit{Fo}(x) - \mathit{Fo}(i + x)}{\mathrm{sinh}(\pi \,p
\,(2\,x + i))}} \,dx={\displaystyle  {{\displaystyle \frac {
i}{2\,p}} \, {\displaystyle \sum _{n=\mathrm{\lfloor}
 - p\rfloor}^{ - 1 + \mathrm{\lfloor}p\rfloor}} \,(-1)^{n}\, \left(  \! 
\mathit{Fo} \left(  \! {\displaystyle \frac {{\displaystyle 
} \,i\,(n + 1 - p)}{2\,p}}  \!  \right)  - \mathit{Fo}
 \left(  \! {\displaystyle \frac {{\displaystyle } \,
i\,(n + 1 + p)}{2\,p}}  \!  \right)  \!  \right)   }{}} 
\]
}
\end{equation}
\begin{equation}
\mapleinline{inert}{2d}{Int((Fo(x)^2+Fo(I+x)^2)/sinh(Pi*p*(2*x+I)),x = -infinity .. infinity)
= 1/2*I*Sum((-1)^n*(Fo(1/2*I*(n+1-p)/p)^2+Fo(1/2*I*(n+1+p)/p)^2),n =
floor(-p) .. -1+floor(p))/p;}{%
\[
{\displaystyle \int _{ - \infty }^{\infty }} {\displaystyle 
\frac {\mathit{Fo}^{2}(x) + \mathit{Fo}^{2}(i + x)}{\mathrm{sinh}
(\pi \,p\,(2\,x + i))}} \,dx={\displaystyle {
{\displaystyle \frac {i}{2\,p}} \, {\displaystyle 
\sum _{n=\mathrm{\lfloor} - p \rfloor}^{ - 1 + \mathrm{\lfloor}p\rfloor}} \,(-1)
^{n}\, \left(  \! \mathit{Fo}^{\displaystyle {2}} \left(  \! {\displaystyle \frac {
{\displaystyle } \,i\,(n + 1 - p)}{2\,p}}  \!  \right) 
 + \mathit{Fo} ^{\displaystyle 2}\left(  \! {\displaystyle \frac {
{\displaystyle } \,i\,(n + 1 + p)}{2\,p}}  \!  \right) 
 \,  \right)  }{}} 
\]
}
\end{equation}

\begin{equation}
\mapleinline{inert}{2d}{Int((Fo(x)*Fe(I+x)-Fe(x)*Fo(I+x))/sinh(Pi*p*(2*x+I)),x = -infinity ..
infinity) =
1/2*I*Sum((-1)^n*(Fo(1/2*I*(n+1-p)/p)*Fe(1/2*I*(n+1+p)/p)-Fo(1/2*I*(n+
1+p)/p)*Fe(1/2*I*(n+1-p)/p)),n = floor(-p) .. -1+floor(p))/p;}{%
\maplemultiline{
{\displaystyle \int _{ - \infty }^{\infty }} {\displaystyle 
\frac {\mathit{Fo}(x)\,\mathit{Fe}(i + x) - \mathit{Fe}(x)\,
\mathit{Fo}(i + x)}{\mathrm{sinh}(\pi \,p\,(2\,x + i))}} \,dx=
{\displaystyle \frac {i}{2\,p}} \left( {\vrule 
height1.83em width0em depth1.83em} \right. \!  \! {\displaystyle 
\sum _{n=\mathrm{\lfloor} - p \rfloor}^{ - 1 + \mathrm{\lfloor}p\rfloor}}(-1)^{n}  \\
\times \left(  \! \mathit{Fo} \left(  \! {\displaystyle 
\frac {{\displaystyle } \,i\,(n + 1 - p)}{2\,p}}  \! 
 \right) \,\mathit{Fe} \left(  \! {\displaystyle \frac {
{\displaystyle } \,i\,(n + 1 + p)}{2\,p}}  \!  \right) 
 - \mathit{Fo} \left(  \! {\displaystyle \frac {{\displaystyle 
} \,i\,(n + 1 + p)}{2\,p}}  \!  \right) \,\mathit{Fe}
 \left(  \! {\displaystyle \frac {{\displaystyle } \,
i\,(n + 1 - p)}{2\,p}}  \!  \right)  \!  \right)  \! \! \left. 
{\vrule height1.83em width0em depth1.83em} \right) }
}
\end{equation}
\normalsize
Using the basic functions
\begin{equation}
\mapleinline{inert}{2d}{F(z) = exp(-b*z);}{%
\[
\mathit{F}(z)=e^{( - b\,z)}
\]
}
\end{equation}
\begin{equation}
\mapleinline{inert}{2d}{Fe(z) = cosh(b*z);}{%
\[
\mathit{Fe}(z)=\mathrm{cosh}(b\,z)
\]
}
\end{equation}
\begin{equation}
\mapleinline{inert}{2d}{Fo(z) = sinh(b*z);}{%
\[
\mathit{Fo}(z)=\mathrm{sinh}(b\,z)
\]
}
\end{equation} 

we find results respectively corresponding to each of the above, as follows
\small
\begin{equation}
\mapleinline{inert}{2d}{Int(exp(-b*v^2)*(sin(b*v)*sinh(2*Pi*p*v)*cos(Pi*p)+cos(b*v)*cosh(2*Pi
*p*v)*sin(Pi*p))/(cos(Pi*p)^2-cosh(2*Pi*p*v)^2),v = 0 .. infinity) 
=1/4*1/p*Sum((-1)^n*exp(1/4*b*((n+1)^2-p^2)/p^2),n = floor(-p) ..
-1+floor(p));}{%
\maplemultiline{
{\displaystyle \int _{0}^{\infty }} {\displaystyle \frac {e^{( - 
b\,v^{2})}\,(\mathrm{sin}(b\,v)\,\mathrm{sinh}(2\,\pi \,p\,v)\,
\mathrm{cos}(\pi \,p) + \mathrm{cos}(b\,v)\,\mathrm{cosh}(2\,\pi 
\,p\,v)\,\mathrm{sin}(\pi \,p))}{\mathrm{cos^{2}}(\pi \,p) - 
\mathrm{cosh^{2}}(2\,\pi \,p\,v)}} \,dv \\
={\displaystyle \frac {1}{4\,p}} \,{\displaystyle  {
{\displaystyle \sum _{n=\mathrm{\lfloor} - p \rfloor}^{ - 1 + \mathrm{
\lfloor}p\rfloor}} \,(-1)^{n}\,e^{(\frac {b\,((n + 1)^{2} - p^{2})}{4\,p
^{2}})}}{}}  }
\label{A11}
}
\end{equation}
\begin{equation}
\mapleinline{inert}{2d}{Int(-sin(b)*sinh(2*Pi*p*v)*cos(Pi*p)*sinh(2*b*v)/(cos(Pi*p)^2-cosh(2*
Pi*p*v)^2)+2*cosh(b*v)^2*cos(b)*cosh(2*Pi*p*v)*sin(Pi*p)/(cos(Pi*p)^2-
cosh(2*Pi*p*v)^2),v = 0 .. infinity) =
1/2*1/p*Sum((-1)^n*cos(1/2*b*(n+1-p)/p)*cos(1/2*b*(n+1+p)/p),n =
floor(-p) .. -1+floor(p));}{%
\maplemultiline{
{\displaystyle \int _{0}^{\infty }} \left( \! - {\displaystyle \frac {
\mathrm{sin}(b)\,\mathrm{sinh}(2\,\pi \,p\,v)\,\mathrm{cos}(\pi 
\,p)\,\mathrm{sinh}(2\,b\,v)}{\mathrm{cos^{2}}(\pi \,p) - 
\mathrm{cosh^{2}}(2\,\pi \,p\,v)}}  + {\displaystyle \frac {2\,
\mathrm{cosh}(b\,v)^{2}\,\mathrm{cos}(b)\,\mathrm{cosh}(2\,\pi \,
p\,v)\,\mathrm{sin}(\pi \,p)}{\mathrm{cos^{2}}(\pi \,p) - 
\mathrm{cosh^{2}}(2\,\pi \,p\,v)}}\!\right) dv  \\
={\displaystyle \frac {1}{2\,p}} \,{\displaystyle {
{\displaystyle \sum _{n=\mathrm{\lfloor} - p \rfloor}^{ - 1 + \mathrm{
\lfloor}p\rfloor}} \,(-1)^{n}\,\mathrm{cos}({\displaystyle \frac {b\,(n
 + 1 - p)}{2\,p}} )\,\mathrm{cos}({\displaystyle \frac {b\,(n + 1
 + p)}{2\,p}} )}{}}  }
}
\end{equation}
\begin{equation}
\mapleinline{inert}{2d}{Int(2*cos(b)*cosh(2*Pi*p*v)*sin(Pi*p)*sinh(b*v)^2/(cos(Pi*p)^2-cosh(2
*Pi*p*v)^2)-sin(b)*sinh(2*Pi*p*v)*cos(Pi*p)*sinh(2*b*v)/(cos(Pi*p)^2-c
osh(2*Pi*p*v)^2),v = 0 .. infinity) =
1/2*1/p*Sum(-(-1)^n*sin(1/2*b*(n+1-p)/p)*sin(1/2*b*(n+1+p)/p),n =
floor(-p) .. -1+\lfloor p \rfloor);}{%
\maplemultiline{
{\displaystyle \int _{0}^{\infty }}\left( \! {\displaystyle \frac {2\,
\mathrm{cos}(b)\,\mathrm{cosh}(2\,\pi \,p\,v)\,\mathrm{sin}(\pi 
\,p)\,\mathrm{sinh}(b\,v)^{2}}{\mathrm{cos^{2}}(\pi \,p) - 
\mathrm{cosh^{2}}(2\,\pi \,p\,v)}}  - {\displaystyle \frac {
\mathrm{sin}(b)\,\mathrm{sinh}(2\,\pi \,p\,v)\,\mathrm{cos}(\pi 
\,p)\,\mathrm{sinh}(2\,b\,v)}{\mathrm{cos^{2}}(\pi \,p) - 
\mathrm{cosh^{2}}(2\,\pi \,p\,v)}} \! \right)\,dv \\
=-{\displaystyle \frac {1}{2\,p}} \,{\displaystyle  {
{\displaystyle \sum _{n=\mathrm{\lfloor} - p \rfloor}^{ - 1 + \mathrm{
\lfloor}p\rfloor}} \,(-1)^{n}\,\mathrm{sin}({\displaystyle \frac {b
\,(n + 1 - p)}{2\,p}} )\,\mathrm{sin}({\displaystyle \frac {b\,(n
 + 1 + p)}{2\,p}} )}{}}  }
}
\end{equation}
\begin{equation}
\mapleinline{inert}{2d}{Int(-sin(b)*sinh(2*Pi*p*v)*cos(Pi*p)*sinh(2*b*v)/(cos(Pi*p)^2-cosh(2*
Pi*p*v)^2)/cosh(b*v)+2*cosh(b*v)*(1+cos(b))*cosh(2*Pi*p*v)*sin(Pi*p)/(
cos(Pi*p)^2-cosh(2*Pi*p*v)^2),v = 0 .. infinity) =
1/2*1/p*Sum((-1)^n*(cos(1/2*b*(n+1-p)/p)+cos(1/2*b*(n+1+p)/p)),n =
floor(-p) .. -1+floor(p));}{%
\maplemultiline{
{\displaystyle \int _{0}^{\infty }}\left( \!  - {\displaystyle \frac {
\mathrm{sin}(b)\,\mathrm{sinh}(2\,\pi \,p\,v)\,\mathrm{cos}(\pi 
\,p)\,\mathrm{sinh}(2\,b\,v)}{(\mathrm{cos^{2}}(\pi \,p) - 
\mathrm{cosh^{2}}(2\,\pi \,p\,v))\,\mathrm{cosh}(b\,v)}}  
\mbox{} + {\displaystyle \frac {2\,\mathrm{cosh}(b\,v)\,(1 + 
\mathrm{cos}(b))\,\mathrm{cosh}(2\,\pi \,p\,v)\,\mathrm{sin}(\pi 
\,p)}{\mathrm{cos^{2}}(\pi \,p) - \mathrm{cosh^{2}(2\,\pi \,p\,v)
}}} \! \right)dv \\
={\displaystyle \frac {1}{p}} \mathrm{cos}({\displaystyle \frac{b}{2}})
{\displaystyle \sum _{n=\mathrm{\lfloor} - p \rfloor}^{ - 1 + \mathrm{
\lfloor}p\rfloor}} \,(-1)^{n}\,\mathrm{cos}({\displaystyle \frac {b\,(n
 + 1)}{2\,p}} )
 }
}
\end{equation}
\begin{equation}
\mapleinline{inert}{2d}{Int(cos(Pi*p)*(-1+cos(b))*sinh(2*Pi*p*v)*sinh(b*v)/(cos(Pi*p)^2-cosh(
2*Pi*p*v)^2)+sin(b)*cosh(b*v)*cosh(2*Pi*p*v)*sin(Pi*p)/(cos(Pi*p)^2-co
sh(2*Pi*p*v)^2),v = 0 .. infinity) =
-1/4*1/p*Sum((-1)^n*(sin(1/2*b*(n+1-p)/p)-sin(1/2*b*(n+1+p)/p)),n =
floor(-p) .. -1+floor(p));}{%
\maplemultiline{
{\displaystyle \int _{0}^{\infty }}\left( \! {\displaystyle \frac {
\mathrm{cos}(\pi \,p)\,( - 1 + \mathrm{cos}(b))\,\mathrm{sinh}(2
\,\pi \,p\,v)\,\mathrm{sinh}(b\,v)}{\mathrm{cos^{2}}(\pi \,p) - 
\mathrm{cosh^{2}}(2\,\pi \,p\,v)}}  + {\displaystyle \frac {
\mathrm{sin}(b)\,\mathrm{cosh}(b\,v)\,\mathrm{cosh}(2\,\pi \,p\,v
)\,\mathrm{sin}(\pi \,p)}{\mathrm{cos^{2}}(\pi \,p) - \mathrm{
cosh^{2}}(2\,\pi \,p\,v)}} \right) \! dv  \\
={\displaystyle \frac {1}{2\,p}} \mathrm{sin}({\displaystyle \frac{b}{2}})
{\displaystyle \sum _{n=\mathrm{\lfloor} - p \rfloor}^{ - 1 + \mathrm{
\lfloor}p\rfloor}} \,(-1)^{n}\,\mathrm{cos}({\displaystyle \frac {b\,(n
 + 1)}{2\,p}} )
 }
}
\end{equation}
\begin{equation}
\mapleinline{inert}{2d}{Int(cos(b)*sin(b)*sinh(2*Pi*p*v)*cos(Pi*p)*sinh(2*b*v)/(cos(Pi*p)^2-c
osh(2*Pi*p*v)^2)-cosh(2*Pi*p*v)*sin(Pi*p)*(cos(b)^2*cosh(2*b*v)-1)/(co
s(Pi*p)^2-cosh(2*Pi*p*v)^2),v = 0 .. infinity) =
1/4*1/p*Sum((-1)^n*(sin(1/2*b*(n+1-p)/p)^2+sin(1/2*b*(n+1+p)/p)^2),n =
floor(-p) .. -1+floor(p));}{%
\maplemultiline{
{\displaystyle \int _{0}^{\infty }} \left( \! {\displaystyle \frac {
\mathrm{cos}(b)\,\mathrm{sin}(b)\,\mathrm{sinh}(2\,\pi \,p\,v)\,
\mathrm{cos}(\pi \,p)\,\mathrm{sinh}(2\,b\,v)}{\mathrm{cos^{2}}(\pi 
\,p) - \mathrm{cosh^{2}}(2\,\pi \,p\,v)}}  
\mbox{} - {\displaystyle \frac {\mathrm{cosh}(2\,\pi \,p\,v)\,
\mathrm{sin}(\pi \,p)\,(\mathrm{cos^{2}}(b)\,\mathrm{cosh}(2\,b\,
v) - 1)}{\mathrm{cos^{2}}(\pi \,p) - \mathrm{cosh^{2}}(2\,\pi \,p\,v)
}} \! \right ) dv \\
={\displaystyle \frac {1}{4\,p}} \,{\displaystyle  {
{\displaystyle \sum _{n=\mathrm{\lfloor} - p \rfloor}^{ - 1 + \mathrm{
\lfloor}p\rfloor}} \,(-1)^{n}\,(\mathrm{sin}^2({\displaystyle \frac {b\,(n
 + 1 - p)}{2\,p}} ) + \mathrm{sin}^2({\displaystyle \frac {b\,(
n + 1 + p)}{2\,p}} ))}{}}  }
}
\end{equation}
\begin{equation}
\mapleinline{inert}{2d}{sin(Pi*p)*Int(cosh(2*Pi*p*v)/(cos(Pi*p)^2-cosh(2*Pi*p*v)^2),v = 0 ..
infinity) = 1/4*Sum((-1)^n,n = floor(-p) .. -1+floor(p))/p;}{%
\[
\mathrm{sin}(\pi \,p)\,{\displaystyle \int _{0}^{\infty }} 
{\displaystyle \frac {\mathrm{cosh}(2\,\pi \,p\,v)}{\mathrm{cos^{2}}(
\pi \,p) - \mathrm{cosh^{2}}(2\,\pi \,p\,v)}} \,dv=
-{\displaystyle \frac{1}{8\,p}\,\left( (-1)^{\lfloor p \rfloor}-(-1)^{\lfloor-p\rfloor} \right)} 
\]
\label{known}
}
\end{equation}
\newline

\normalsize
Except for \eqref{A11} which reproduces \eqref{C10b}, and \eqref{known}, which corresponds to \eqref{X1}, these results, although not independent, are believed to be new, and clearly reduce to more useful results in special cases of the parameters $b$ and $p$, whose range is initially always chosen such that the integral representations exist. This is left as an exercise for the reader. A small number of related tabular results can be found in \cite [Eq. 3.543(2) and Sections 3.983 and 3.984] {GR}. Except for \eqref{A11}, Maple \cite{Maple} is able to analytically evaluate the exact sums given for the right-hand sides of each of these results. But the forms given here are generally more succinct.

\end{flushleft}

\noindent

\begin{thebibliography}{100}
\bibliographystyle{ieeetr}
\bibitem{Remarkable}
M.L. Glasser,
{\it A Remarkable Definite Integral}, 
{\url{http://arxiv.org/pdf/1308.6361v2.pdf}}

\bibitem{Berndt}
Bruce C. Berndt,
{\it Ramanujan's Notebooks Part IV}, 
[Springer Verlag, New York(1993)], p.317

\bibitem{Tichmarsh}
E.C. Tichmarsh,
{\it The Riemann Zeta Function}, 
[Cambridge Univ. Press, London 1930]

\bibitem{Hardy}
G.H. Hardy, 
{\it Ramanujan},
[Chelsea Publishing Co. New York, 1940] Chap. XI

\bibitem{GRHO}
W. Gr{\"o}bner, N. Hofreiter,
{\it IntegralTafel, (Bestimmte Integrale)},
Springer-Verlag, Vienna and Innsbruck, (1958), Volume 2.

\bibitem{Maple}
© Maplesoft, a division of Waterloo Maple Inc. 2014.
\url{http://www.maplesoft.com/products/maple/}

\bibitem{Math}
Wolfram Research, Inc., Mathematica, Version 9.0, Champaign, IL (2012).

\bibitem{Luke}
Y. Luke
{\it The Special Functions and their Approximations, Vol. 1},
Academic press, (1969).


\bibitem{Bateman}
H. Bateman
{\it Higher Transcendental Functions Volume 1}
McGraw-Hill (1953).


\bibitem{GR}
I. Gradshteyn and I. Ryzhik,
{\it Tables of Integrals, Series and Products},
Academic Press (1980) Corrected and Enlarged Edition.

\bibitem{G&M}
M.L. Glasser and M. Milgram,
{\it On Quadratic Gauss Sums and Variations Thereof},
submitted for publication, (2014).

\end{thebibliography}
\end{document}